\documentclass[a4paper,12pt,oneside]{article}
\usepackage{custom-style}

\title{%
Random sketching of operators with application to learning preconditioners
}
\author{%
  O. Balabanov${}^{1}$,
  A. Nouy${}^{2}$,
  A. Pasco${}^{2}$}
\date{\medskip%
  \small %
  ${}^1$  International Computer Science Institute, University of California, Berkeley \\
  Lawrence Berkeley National Laboratory \\ 
  \texttt{olegbalabanov@gmail.com} \\
  ${}^2$ Centrale Nantes, Nantes Universit\'e, \\ Laboratoire de Mathématiques Jean Leray UMR CNRS 6629 \\
  \texttt{anthony.nouy@ec-nantes.fr}\\
  \texttt{pasco.alexandre@proton.me}
}

\begin{document}
\maketitle

\begin{abstract}

We propose a new random sketching approach for embedding high-dimensional Hilbert-Schmidt operators, using random input-output pairs.
Such operator can then be approximated in a low-dimensional subspace of operators by solving a small least-squares problem.
To achieve computational efficiency, we introduce a structured
random map, composed of three random matrices.
We provide rigorous conditions under which subspaces of operators are accurately embedded with high probability.
The framework is flexible, as the random matrices may be adapted to the operator structure and the computational environment.

As an application, we consider the construction of preconditioners for high-dimensional linear equations.
We derive a rigorous characterization of preconditioner quality through the discrepancy between the preconditioned operator and an optimal baseline, which can be tailored to a linear approximation space for the solution.
We show that this quantity can be efficiently minimized within the proposed framework, especially for parameter separable linear equations.
We then establish rigorous high-probability bounds on the quasi-optimality error of the preconditioned Galerkin projection and on the accuracy of a preconditioned residual-based error estimator when the sketch dimensions are sufficiently large.
Numerical experiments on an acoustic wave scattering benchmark demonstrate the effectiveness of the method.

\end{abstract}





\section{Introduction}
\label{sec:introduction preconditioners}

Randomization and dimensionality reduction play a central role in scientific computing and machine learning. 
Modern randomized methods make effective use of computational resources and can tackle problems at scales far beyond the reach of deterministic techniques, with strong guarantees of accuracy \cite{martinssonRandomizedNumericalLinear2020}.

We propose a matrix-free operator learning approach based on random sketching.
It uses efficient linear random embeddings $\bfTheta$ into low-dimensional Euclidean spaces $\bbK^k$ (with $\bbK = \bbR$ or $\bbC$)  that, with high probability, approximately preserve inner products.
More precisely, for all pairs of elements $\bfX$ and $\bfY$ from subspaces of interest we have
\[
    \innerp{\bfX}{\bfY} \approx \innerp{\bfTheta(\bfX)}{\bfTheta(\bfY)}.
\]
This property then allows us to perform the original computationally intensive task, such as least-squares fitting, in $\bbK^k$ at reduced computational cost.
We refer to \cite{woodruffSketchingToolNumerical2014} for a broad discussion on subspace embeddings.
Classical sketching typically considers $\bfX$ and $\bfY$ as Euclidean vectors endowed with the canonical inner product. 
In this work, we extend this framework and introduce structured embeddings $\bfTheta$ for high-dimensional Hilbert-Schmidt operators.

Designing efficient randomized methods suited to matrix-free settings is crucial in large scale problems.
A few related examples are the following.
In \cite{hutchinsonStochasticEstimatorTrace1989,bujanovicNormTraceEstimation2021,perssonImprovedVariantsHutch2022,zvonekContHutchStochasticTrace2025} randomized trace or norm estimators were considered, although subspace embeddings were not incorporated.
In \cite{elsworthConversionsBarycentricRKFUN2019,guttelRandomizedSketchingNonlinear2024,chenRobustRandomizedIndicator2024,bujanovicSubspaceEmbeddingRandom2025} nonlinear eigenvalue problems were considered.
In particular, \cite{bujanovicSubspaceEmbeddingRandom2025} derives subspace embedding properties for Khatri-Rao product of Gaussian embeddings, closely related to our approach.
Finally, \cite{zahmInterpolationInverseOperators2016} considers preconditioning of parameter-dependent equations, where parameter-dependent preconditioners are constructed from a linear subspace by approximately minimizing an error in Frobenius norm using a randomized approach also closely related to ours.

In this work, we apply our operator sketching method to build preconditioners for approximating the solution $\bfu$ in some Hilbert space $U$, $\dim U \gg 1$, of a large-scale parameter-dependent linear equation
\begin{equation}
\label{equ:def fom}
    \bfA \bfu = \bfb,
\end{equation}
with a linear operator $\bfA$ and right-hand side $\bfb$.
In practice, the space $U$ often arises from the discretization of some function space, using for example the finite elements or finite volumes methods.
Given a subspace $U_r\subset U$, we address the problem of efficiently computing an approximation of $\bfu$ in this subspace.
This type of problem typically arises in linear \emph{model order reduction} (MOR) methods for parameterized PDEs, see for example \cite{haasdonkChapter2Reduced2017}, which we investigate in the present work, or in Krylov subspace methods.
A classical method for computing such an approximation is the Galerkin projection.
However, when $\bfA$ is ill-conditioned, this projection can be far from the orthogonal projection, and residual-based error estimators may provide very crude estimates of the error in $U$-norm.
This may result in large approximation error even when $\bfu$ can be well approximated in $U_r$, and thus constitutes one of the central issues in MOR and Krylov methods.
Hence, our goal is to construct a preconditioner $\bfP$ to address these two problems.
Inspired by \cite{zahmInterpolationInverseOperators2016}, we investigate a linear approximation of the inverse of $\bfA$, so that $\bfP$ is taken as the solution to
\begin{equation}
\label{equ:proj error A onto precond space}
    \min_{\bfP \in \spanv{\bfY_1, \cdots, \bfY_p}}
    \|\bfI - \bfP \bfA\|_{HS},
\end{equation}
for some operators $\bfY_i$ and some Hilbert-Schmidt norm $\|\cdot\|_{HS}$, which is a linear least-squares problem with $p$ degrees of freedom.
In the parameter dependent setting, where both $\bfA$ and $\bfb$ are parameter-dependent, we construct the approximation space for the preconditioner using a classical greedy algorithm, where $\bfY_i$ is taken as $\bfA(\xi_i)^{-1}$ for some parameter value $\xi_i$ at which an error criterion, based on the aforementioned error measure, is maximal.
Other strategies for selecting the interpolation points can be found in \cite{zahmInterpolationInverseOperators2016}.

Constructing preconditioners in similar settings was considered in \cite{elmanPreconditioningTechniquesReduced2015,zahmInterpolationInverseOperators2016,pasettoReducedOrderModelbased2017,santoMultiSpaceReduced2018,lindsayPreconditionedLeastsquaresPetrov2022}.
In \cite{pasettoReducedOrderModelbased2017,santoMultiSpaceReduced2018} the reduced space was used to construct preconditioners for accelerating iterative methods used to solve \eqref{equ:def fom}.
In \cite{elmanPreconditioningTechniquesReduced2015} preconditioners were constructed for accelerating iterative methods used to solve the reduced Galerkin system.
In particular, \cite{zahmInterpolationInverseOperators2016} constructs preconditioners to improve the projection onto the reduced space, improve the accuracy of residual-based error estimators, and reuse factorized operators computed during the construction of the reduced space.


\subsection{Contributions}

There are two major contributions in this work.
The first major contribution is a new random sketching method for Hilbert-Schmidt (HS) operators between finite-dimensional spaces, presented in \Cref{sec:random sketching for linear operators}.
We introduce a random embedding for matrices,
\[
    \bfTheta(\bfX)
    := \bfGamma \mathrm{vec}(\bfOmega \bfX \bfSigma^*),
\]
where $\bfOmega$, $\bfSigma$ and $\bfGamma$ are classical random matrices such as Gaussian, structured or sparse.
We show that if $\bfOmega$, $\bfSigma$ and $\bfGamma$ are oblivious subspace embeddings for the canonical inner product, then $\bfTheta$ is an oblivious subspace embedding for the Frobenius inner product.
An important property of this embedding is that it can be applied to matrices given in matrix-free form, meaning that it only requires matrix-vector products with $\bfX$ to be computed.
More precisely, we show that any element of a $d$-dimensional subspace of matrices can be embedded with distortion $\varepsilon$ and probability $1-\delta$ by applying it to only $\calO(\varepsilon^{-2} (d + \log(1/\delta)))$ independent random vectors.
We also provide non-asymptotic bounds on the sketch sizes.
Then, in \Cref{subsec:ose in hs spaces}, we inspire from \cite{balabanovRandomizedLinearAlgebra2019} to extend this embedding to HS operators on finite dimensional Hilbert spaces with arbitrary inner products.

The second major contribution of this work is presented in \Cref{sec:measures of quality of a preconditioner}.
In \Cref{subsec:general purpose preconditioner} we extend the approach from \cite{zahmInterpolationInverseOperators2016}, which is used to construct a general purpose preconditioner, to HS operators.
The important property of the discrepancy measure $\|\bfI - \bfP \bfA\|$ in HS norm is that it is an upper bound of the discrepancy measure in operator norm used in theoretical results, while being well suited to random sketching approaches, as described in \Cref{subsec:sketched measures of quality}.
The problem is that it can overestimate the operator norm by a large factor, up to $n^{1/2}$ with $n := \mathrm{dim}(U)$.

This problem can be circumvented when constructing preconditioners tailored to a reduced space, which is denoted as a MOR purpose, as we propose in \Cref{subsec:mor-purpose preconditioner}.
We introduce new discrepancy measures $\|\bfI - \bfP \bfA\|_{U,U_r}$ and $\|\bfI - \bfP \bfA\|_{U_r,U_r}$ associated to an $r$-dimensional subspace $U_r$, $r\ll n$, where $\|\cdot\|_{U,U_r}$ and $\|\cdot\|_{U_r,U_r}$ are operator seminorms introduced in \Cref{sec:norms on linear operators}.
We show that the preconditioned Galerkin projection $\bfu_r$ satisfies
\[
    \|\bfu - \bfu_r\|_U \leq 
    (1 + \frac{\|\bfI - \bfP \bfA\|_{U, U_r}}{1 - \|\bfI - \bfP \bfA\|_{U_r, U_r}})
    \|\bfu - \bfPi_{U_r} \bfu\|_U,
\]
assuming that $\|\bfI - \bfP \bfA\|_{U_r, U_r} <1$.
We also introduce a new preconditioned residual-based error estimator which makes use of an $m$-dimensional space $U_m$ with $r \leq m \ll n$, which contains $U_r$.
We show that under some assumptions, we can quantify the accuracy of this estimator using $\|\bfI - \bfP \bfA\|_{U,U_m}$ and $\|\bfI - \bfP \bfA\|_{U_m,U_m}$.
An important property of these new discrepancy measures is that their Hilbert-Schmidt counterparts can overestimate them by only a small factor, up to $r^{1/2}$ or $m^{1/2}$.
These Hilbert-Schmidt counterparts are also well suited for random sketching approaches, as described in \Cref{subsec:sketched measures of quality}.

Then, in the parameter-dependent setting with parameter separable operator, although solving \eqref{equ:proj error A onto precond space} using the normal equation can be done online efficiently, this method can suffer from round-off errors, as pointed out in \cite{buhrNumericallyStablePosteriori2014,casenaveAccurateOnlineefficientEvaluation2014,balabanovRandomizedLinearAlgebra2019}.
Also, the corresponding offline costs may be prohibitive.
These problems are circumvented by the new sketched discrepancy measures introduced in \Cref{subsec:sketched measures of quality}, which are based on our operator sketching approach from \Cref{sec:random sketching for linear operators}.
Indeed, we approximate the solution of \eqref{equ:proj error A onto precond space} by replacing the Hilbert-Schmidt norm by its sketched counterpart, resulting in a small sketched least-squares problem, which can be solved efficiently online using stable methods, and whose offline cost is reasonable.
We detail in \Cref{sec:practical aspects preconditioners} the practical aspects and computational costs associated to our approach in the case where $\bfA$ is parameter separable with $m_{\bfA}$ affine terms.
We show that our approach is online efficient, as constructing $\bfP(\xi)$ online is done by assembling and solving a small linear least-squares problem of size $k \times p$, which costs in total $\calO(kp(p+m_{\bfA}))$ flops.
We also show that, using structured embeddings and parameter independent $\bfY_i$, the offline costs associated to the sketched counterpart of \eqref{equ:proj error A onto precond space} for our general purpose and our MOR purpose discrepancy measures are respectively
\[
    \calO\left(kp(T_{\bfY} + m_{\bfA} n\log(n))\right)
    \quad
    \text{and}
    \quad
    \calO\left(\min(k,m) p (T_{\bfY} + \min(k,m) m_{\bfA} n\log(n)) \right),
\]  
where $T_{\bfY}$ denotes the cost of applying $\bfY_i$ to a vector.
We illustrate our approach on an acoustic wave scattering numerical example in \Cref{sec:numerics preconditioners}.
Implementation of our approach is freely available at \url{https://github.com/alexandre-pasco/rla4mor}.

A preliminary version of this work appeared as a preprint \cite{balabanovPreconditionersModelOrder2021}.
The current article significantly extends this earlier version, in particular in the theoretical analysis of random sketching of operators, the construction of preconditioners with control in Hilbert Schmidt norms, and numerical evaluation.

\subsection{Outline}

The paper is organized as follows.
In section \Cref{sec:norms on linear operators} we introduce some notations as well as norms and seminorms of linear operators between finite dimensional Hilbert spaces.
In \Cref{sec:random sketching for linear operators} we introduce new random embeddings for such operators, and we discuss on other existing methods.
In \Cref{sec:measures of quality of a preconditioner} we introduce measures of quality of a preconditioner, depending on the purpose, as well as their sketched versions.
In \Cref{sec:practical aspects preconditioners} we provide practical details for the parameter-dependent setting with parameter separable operator and right-hand side.
In \Cref{sec:numerics preconditioners} we illustrate the described methods with an acoustic wave scattering numerical example.
Finally, in \Cref{sec:conclusion preconditioners} we summarize the analysis and observations, and we discuss perspectives.


\section{Preliminary on norms of linear operators}
\label{sec:norms on linear operators}

In this section, we introduce norms and seminorms of linear operators between finite dimensional Hilbert spaces over the field $\bbK$ ($\bbR$ or $\bbC$), identified with matrices.
We let $V = (\bbK^{\dim V}, \innerp{\cdot}{\cdot}_V)$ be a finite-dimensional Hilbert spaces equipped with the inner product $\innerp{\cdot}{\cdot}_V$ and associated norms $\|\cdot\|_V$.
Let $V'$ be the dual of $V$ and let $\bfR_V : V \rightarrow V'$ be the self-adjoint positive definite matrix such that $\innerp{\bfv_1}{\bfv_2}_V = \innerp{\bfR_V \bfv_1}{\bfv_2} = \bfv_1^* \bfR_V \bfv_2$ for all $\bfv_1,\bfv_2\in V$, where $\innerp{\cdot}{\cdot}$ denotes the duality bracket between $V'$ and $V$.
We write $V \equiv \bbK^{\dim V}$ when equipped with the canonical inner product $\innerp{\cdot}{\cdot}_2$ and associated inner product $\|\cdot\|_2$, in which case $V'$ is identified with $V$.
We define in the same way $W = (\bbK^{\dim W}, \innerp{\cdot}{\cdot}_W)$ where $\innerp{\cdot}{\cdot}_W = \innerp{\bfR_W\cdot}{\cdot}$.

In the present work, we focus on two norms.
Let $\bfG : V \rightarrow W$ be a linear operator. The \emph{operator norm} of $\bfG$, denoted by $\|\bfG\|_{V, W}$, is defined by
\begin{equation}
\label{equ:def operator norm}
    \|\bfG\|_{V, W} 
    := \max_{\bfv \in V \setminus \{0\}}
    \frac{\|\bfG \bfv\|_W}{\|\bfv\|_V},
\end{equation}
and is equal to its largest singular value.
The second norm we focus on is the \emph{Hilbert-Schmidt} (HS) norm of $\bfG$, denoted by $\|\bfG\|_{HS(V, W)}$, which is induced by the inner product $\innerp{\cdot}{\cdot}_{HS(V,W)}$, defined by
\begin{equation}
\label{equ:def hs norm}
    \|\bfG\|_{HS(V, W)}
    := \innerp{\bfG}{\bfG}_{HS(V, W)}^{1/2},
    \quad 
    \innerp{\bfG_1}{\bfG_2}_{HS(V,W)}
    := \sum_{i=1}^{\dim V}
    \innerp{\bfG_1 \bfv_i}{ \bfG_2 \bfv_i}_W,
\end{equation}
for any linear operators $\bfG_1, \bfG_2 : V \rightarrow W$ and orthonormal basis $(\bfv_i)_{1\leq i\leq \dim V}$ of $V$.
The HS norm and inner product can be expressed with the classical Frobenius norm $\|\cdot\|_F$ and associated inner product $\innerp{\cdot}{\cdot}_F$, 
\begin{equation}
\label{equ:hs as frobenius}
\begin{aligned}
    & \|\bfG\|_{HS(V, W)} 
    = \|\bfQ_W \bfG \bfR_V^{-1} \bfQ_V^*\|_{F}, \\
    & \innerp{\bfG_1}{\bfG_2}_{HS(V,W)}
    = \innerp{\bfR_W \bfG_1}{\bfG_2 \bfR_V^{-1}}_{F}
    = \Trace{\bfG_1^* \bfR_W \bfG_2 \bfR_V^{-1}},
\end{aligned}
\end{equation}
for any matrices $\bfQ_V \in \bbK^{p \times \dim V}$ and $\bfQ_W \in \bbK^{q \times \dim W}$ such that $\bfQ_V^* \bfQ_V = \bfR_V$ and $\bfQ_W^* \bfQ_W = \bfR_W$.
Such matrices can be computed for example by Cholesky decompositions of $\bfR_V$ and $\bfR_W$, or a more efficient domain decomposition approach \cite{balabanovRandomizedLinearAlgebra2019}.
Note that it is convenient to see these matrices as operators $\bfQ_V : V \rightarrow \bbK^p$ and $\bfQ_W : W \rightarrow \bbK^q$.
We define the space $HS(V,W) := (\bbK^{\dim W \times \dim V}, \innerp{\cdot}{\cdot}_{HS(V,W)})$.
Noting that $\|\cdot\|_{HS(\bbK^{\dim V},\bbK^{\dim W})} = \|\cdot\|_F$, we write $HS(\bbK^{\dim V}, \bbK^{\dim W}) \equiv \bbK^{\dim W \times \dim V}$.

The HS norm defines an inner product on the space of linear operators from $V$ to $W$, thus solving \eqref{equ:proj error A onto precond space} is equivalent to solving a linear least-squares problem in $\bbK^{\dim V \times \dim W}$ with $p$ degrees of freedom.
Also, in the finite dimensional setting, it is equivalent to the operator norm, as stated in \Cref{prop:hs norm equivalent op norms}.

\begin{proposition}
\label{prop:hs norm equivalent op norms}
    Let $\bfG : V \mapsto W$ a linear operator, with $V$ and $W$ finite dimensional Hilbert spaces.
    It holds
    \begin{equation}
    \label{equ:hs norm equivalent op norms}
        \frac{1}{\sqrt{\min(\dim V, \dim W)}} \|\bfG\|_{HS(V, W)} 
        \leq \|\bfG\|_{V, W}
        \leq \|\bfG\|_{HS(V, W)},
    \end{equation}
    with $\|\cdot\|_{HS(V, W)}$ and $\|\cdot\|_{V, W}$ defined in \eqref{equ:def hs norm} and \eqref{equ:def operator norm} respectively.
\end{proposition}
\begin{proof}
    Let us assume that $\dim V \leq \dim W$.
    This does not cause any loss of generality, as we can apply the analysis to $\bfG^* : W' \rightarrow V'$ to handle the case $ \dim W < \dim V$.
    Indeed, from \eqref{equ:hs as frobenius} we have that $\|\bfG^*\|_{HS(W', V')} = \|\bfG\|_{HS(V,W)}$ and from \eqref{equ:def operator norm} we have that $\|\bfG^*\|_{W', V'} = \|\bfG\|_{V,W}$.
    Moreover, $\dim W' = \dim W$ and $\dim V' = \dim V$.
    
    Now, consider orthonormal bases $(\bfv_i)_{1\leq i\leq \dim V}$ and $(\bfw_i)_{1\leq i\leq \dim W}$ of respectively $V$ and $W$.
    Firstly, let us show the right inequality in \eqref{equ:hs norm equivalent op norms}.
    Let $\bfv \in V$ with $\bfv = \sum_i v_i \bfv_i$ and $\bfG \bfv_i = \sum_j g_{ij} \bfw_j \in W$.
    By expanding $\|\bfG \bfv\|_W^2$, using Cauchy-Schwarz inequality, factorizing by $\|\bfv\|_V^2$, noticing that $\sum_j g_{ij}^2 = \|\bfG \bfv_i\|_W^2$ and using the definition \eqref{equ:def hs norm}, we obtain 
    \[
    \begin{aligned}
        \|\bfG \bfv\|_W^2
        &= \|\sum_i v_i \bfG \bfv_i\|_W^2
        = \|\sum_i v_i \sum_j g_{ij}\bfw_j\|_W^2
        = \|\sum_j (\sum_i v_i g_{ij}) \bfw_j\|_W^2
        = \sum_j (\sum_i v_i g_{ij})^2
        \\ 
        & \leq \sum_j (\sum_i v_i^2) (\sum_i g_{ij}^2)
        = \|\bfv\|_V^2 \sum_{i, j} g_{ij}^2
        = \|\bfv\|_V^2 \sum_j \|\bfG \bfv_i\|_W^2
        = \|\bfv\|_V^2 \|\bfG\|_{HS(V,W)}.
    \end{aligned}
    \]
    Then the definition of the operator norm in \eqref{equ:def operator norm} yields the desired right inequality in \eqref{equ:hs norm equivalent op norms}.
    Secondly, let us show the left inequality in \eqref{equ:hs norm equivalent op norms}.
    Using the definition \eqref{equ:def hs norm}, the definition \eqref{equ:def operator norm} and the fact that $\|\bfv_i\|_V = 1$ for all $1\leq i\leq \dim V$, we obtain
    \[
        \|\bfG\|_{HS(V, W)}^2
        = \sum_{i=1}^{\dim V} \|\bfG \bfv_i\|_W^2
        \leq \sum_{i=1}^{\dim V} 
        \|\bfG\|_{V, W}^2 \|\bfv_i\|_V^2
        = \|\bfG\|_{V, W}^2 \dim V.
    \]
    Then dividing by $\dim V >0$ and taking the square root yields the desired left inequality in \eqref{equ:hs norm equivalent op norms}.
\end{proof}

A potential problem of the HS norm as a surrogate to the operator norm is that the left inequality in \eqref{equ:hs norm equivalent op norms} may be very coarse when $\dim V = \dim W \gg 1$.
Finally, for any subspaces $\tilde V \subset V$ and $\tilde W \subset W$, we define the seminorms
\begin{equation}
\label{equ:def operator hs seminorm}
    \|\bfG\|_{\tilde V, \tilde W} := \|\bfPi_{\tilde W} \bfG \bfPi_{\tilde V}\|_{V, W},
    \quad
    \|\bfG\|_{HS(\tilde V, \tilde W)} := \|\bfPi_{\tilde W} \bfG \bfPi_{\tilde V}\|_{HS(V, W)},
\end{equation}
where $\bfPi_{\tilde V} : V \rightarrow \tilde V$ and $\bfPi_{\tilde W} : W \rightarrow \tilde W$ denote the orthogonal projectors on $\tilde V$ and $\tilde W$ respectively.
It is important to note that those seminorms also satisfy \Cref{prop:hs norm equivalent op norms}.
It is also worth noting that $\|\cdot\|_{HS(\tilde V, \tilde W)}$ defines a corresponding semi-inner-product, hence \eqref{equ:proj error A onto precond space} with this seminorm is a linear least-squares problem.


\section{Random sketching for Hilbert-Schmidt operators}
\label{sec:random sketching for linear operators}

In this section, we introduce new random embeddings for Hilbert--Schmidt operators between finite-dimensional Hilbert spaces. 
These embeddings have two key properties. 
Firstly, they can be applied using only applications of the operators to vectors. 
Secondly, they are defined in a flexible form using standard subspace embeddings on $\bbK^n$ and are shown to inherit their oblivious subspace embedding property. 
As a result, they can be combined with any classical random embedding matrices. 

The outline of this section is as follows.
In \Cref{subsec:ose} we recall the concept of oblivious subspace embedding.
In \Cref{subsec:ose in frobenius space} we introduce our new random embedding on spaces of matrices, which we generalize to HS operators in \Cref{subsec:ose in hs spaces}.
Finally, in \Cref{subsec:comparison with existing methods} we analyze the differences between our approach and other existing approaches.

\subsection{Oblivious subspace embeddings}
\label{subsec:ose}

For $Z$ a finite dimensional vector space equipped with an inner-product $\innerp{\cdot}{\cdot}_Z$ and a linear map $\bfTheta : Z \rightarrow \bbK^k $ with $k \leq \dim Z$, we define the following semi-inner product and seminorm on $Z$,
\begin{equation}
\label{equ:def sketched norm}
    \innerp{\bfx}{\bfy}^{\bfTheta}_Z
    := \innerp{\bfTheta(\bfx)}{\bfTheta(\bfy)}_2,
    \quad
    \|\bfx\|_Z^{\bfTheta} := \|\bfTheta(\bfx)\|_2.
\end{equation}
The idea here is to approximate the norms and inner products of high dimensional vectors in $Z$ by the norms and inner products of their low dimensional images by $\bfTheta$.
For this approximation to be accurate, we will ask $\bfTheta$ to satisfy some quasi-isometry property.
Obviously, if $k < \dim Z$, we cannot ask for quasi-isometry over $Z$ entirely.
However, we can ask it over some low-dimensional subspace $\tilde Z \subset Z$ with $\dim \tilde Z \leq k$, which leads to the definition of \emph{subspace embedding} in \Cref{def:subspace embedding}.

\begin{definition}
\label{def:subspace embedding}
    A linear map $\bfTheta : Z \rightarrow \bbK^k$ is called an $\varepsilon$-embedding for a subspace $\tilde Z \subset Z$ if 
    \begin{equation}
        \forall \bfx,\bfy \in \tilde Z, ~
            |\innerp{\bfx}{\bfy}_Z 
            - \innerp{\bfx}{\bfy}_Z^{\bfTheta}|
        \leq \varepsilon \|\bfx\|_Z \|\bfy\|_Z.
    \end{equation}
\end{definition}

It is rather straightforward to build such an embedding for a fixed subspace $\tilde Z$, as one can simply take $\bfx \mapsto (\innerp{\bfz_i}{\bfx}_Z)_{1\leq i\leq \dim \tilde Z}$, with an orthonormal basis $(\bfz_i)_{1\leq i\leq \dim \tilde Z}$ of $\tilde Z$, which is actually an isometry from $\tilde Z$ to $\bbK^{\dim(\tilde Z)}$, thus it is a $0$-embedding for $\tilde Z$.
However, there are situations where we want $\bfTheta$ to be an $\varepsilon$-embedding simultaneously for an entire family of low-dimensional subspaces, which may not even be known before constructing $\bfTheta$.

A classical way to build such embedding is to take $\bfTheta$ as a realization of some random linear map, whose distribution shall be carefully chosen to ensure theoretical or practical properties.  
For a fixed subspace $\tilde Z \subset Z$, we will then ask for $\bfTheta$ to be a subspace embedding with high probability.
This leads to the definition of a \emph{data-oblivious subspace embedding}, or \emph{oblivious subspace embedding}, in \Cref{def:data oblivious subspace embedding}.

\begin{definition}
\label{def:data oblivious subspace embedding}
    A random linear map $\bfTheta : Z \rightarrow \bbK^k$ is called a $(\varepsilon, \delta, d)$ oblivious $Z \rightarrow \bbK^k$ subspace embedding if for any subspace $Z_d \subset Z$ with $\dim Z_d = d$, it holds 
    \begin{equation}
    \label{equ:data oblivious subspace embedding}
        \Proba[\bfTheta]{
            \forall \bfx,\bfy \in Z_d, ~
            |\innerp{\bfx}{\bfy}_Z 
            - \innerp{\bfx}{\bfy}_Z^{\bfTheta}|
            \leq \varepsilon \|\bfx\|_Z \|\bfy\|_Z
        } \geq 1 - \delta.
    \end{equation}
\end{definition}

It remains now to know what kind of random linear maps actually yield oblivious embeddings.
In the case where $Z = \bbR^n$ we know from \cite{woodruffSketchingToolNumerical2014} several distributions of random matrices $\bfTheta \in \bbR^{k \times n}$ that satisfy \eqref{def:data oblivious subspace embedding}.
It is for example the case for the rescaled Gaussian and Rademacher distributions with $k = \calO(\varepsilon^{-2} (d + \log(1/\delta)))$, for the partial subsampled randomized Hadamard transform (P-SRHT) with $k = \calO(\varepsilon^{-2} \log(d\delta^{-1}) (d + \log(\varepsilon^{-1}\delta^{-1})))$, see \cite{cohenOptimalApproximateMatrix2015}, or for sparse embeddings with $k=\calO(\varepsilon^{-2} (d + \log(\varepsilon^{-1}\delta^{-1})))$, see \cite{chenakkodOptimalObliviousSubspace2025}.
A rescaled Gaussian matrix of size $k\times n$ has Gaussian i.i.d. entries with mean 0 and variance $k^{-1}$. 
A P-SRHT matrix of size $k\times n$ is defined as the first $n$ columns of $k^{-1/2}(\bfR \bfH_s \bfD)$, where $s$ is the power of $2$ such that $n \leq s < 2n$, $\bfR \in \bbR^{k\times s}$ are the first $k$ rows of a uniform random permutation of rows of the identity matrix, $\bfH_s \in \bbR^{s\times s}$ is a Walsh-Hadamard matrix, and $\bfD \in \bbR^{s \times s}$ is a diagonal matrix whose diagonal entries are i.i.d random sign flips.
Other classical examples can be found in \cite{woodruffSketchingToolNumerical2014}, while more recent ones were introduced for example in \cite{balabanovBlockSubsampledRandomized2023,bujanovicSubspaceEmbeddingRandom2025}.
Note that explicit non-asymptotic lower bounds on $k$ in this setting have been proven in \cite{balabanovRandomizedLinearAlgebra2019}.

Those results have then been extended in \cite{balabanovRandomizedLinearAlgebra2019} to the case $Z = (\bbK^n, \innerp{\cdot}{\bfR_Z \cdot})$ with $\bfR_Z \in \bbK^{n\times n}$ a self-adjoint positive definite matrix.
Denoting $\bfQ_Z \in \bbK^{s\times n}$ any matrix such that $\bfQ_Z^* \bfQ_Z = \bfR_Z$, if $\bfOmega \in \bbK^{k \times s}$ is a $(\varepsilon, \delta, d)$ oblivious $\bbK^n \rightarrow \bbK^k$ subspace embedding, then $\bfTheta = \bfOmega \bfQ_Z$ is a $(\varepsilon, \delta, d)$ oblivious $Z \rightarrow \bbK^k$ subspace embedding.

\subsection{Embeddings in matrix space equipped with Frobenius norm}
\label{subsec:ose in frobenius space}

Let us consider $Z = \bbK^{q\times p}$ equipped with the Frobenius inner product as introduced in \Cref{sec:norms on linear operators}.
A first idea to embed elements of $Z$ could be to identify this space with $\bbK^{qp}$, which is possible as long as we have explicit access to the elements of $Z$.
However, we may need to embed matrices defined only implicitly, meaning that only matrix-vector multiplication is available.
This is for example the case when embedding inverses or factorizations of large (sparse) matrices.
We thus need to define an embedding $\bfUpsilon$ such that $\bfUpsilon(\bfX)$ can be computed by involving only a moderate number of matrix-vector multiplications of $\bfX\in\bbK^{q\times p}$.
We propose to consider the following sketch for matrices
\begin{equation}
\label{equ:def upsilon}
    \bfUpsilon(\bfX) 
    := \bfGamma \mathrm{vec}(\bfOmega \bfX \bfSigma^*) \in \bbK^{k},
\end{equation}
for some classical embeddings for vectors $\bfSigma \in \bbK^{k_{\bfSigma} \times p}$, $\bfOmega \in \bbK^{k_{\bfOmega} \times q}$ and $\bfGamma \in \bbK^{k \times k_{\bfSigma} k_{\bfOmega}}$.
Note that the main matrix compression of this sketch is performed by $\bfOmega$ and $\bfSigma$, resulting in $\bfOmega \bfX \bfSigma^* \in \bbK^{k_\bfOmega \times k_{\bfSigma}}$, while the role of $\bfGamma$ is to recompress the rather high dimensional vector $\mathrm{vec}(\bfOmega \bfX \bfSigma^*)$. 
An important property of $\bfUpsilon$ is that it inherits oblivious subspace embedding properties of $\bfSigma$, $\bfOmega$ and $\bfGamma$, as we show in the next sections.

\hspace{1em}
\subsubsection{General oblivious subspace embeddings}

\Cref{prop:matrix embedding} shows how $\bfUpsilon$ can inherit oblivious subspace embedding properties from embeddings $\bfSigma$, $\bfOmega$ and $\bfGamma$.

\begin{proposition}
\label{prop:matrix embedding}
    Let $\bfSigma \in \bbK^{k_{\bfSigma}\times p}$ be a $(\varepsilon/4, \delta/(3q), d)$ oblivious $\bbK^p \rightarrow \bbK^{k_{\bfSigma}}$ subspace embedding.
    Let $\bfOmega \in \bbK^{k_{\bfOmega}\times q}$ be a $(\varepsilon/4, \delta/(3k_{\bfSigma}), d)$ oblivious $\bbK^q \rightarrow \bbK^{k_{\bfOmega}}$ subspace embedding.
    Let $\bfGamma \in \bbK^{k\times k_{\bfSigma} k_{\bfOmega}}$ be a $(\varepsilon/4, \delta/3, d)$ oblivious $\bbK^{k_{\bfSigma} k_{\bfOmega}} \rightarrow \bbK^{k}$ subspace embedding.
    Assume that $\bfSigma$, $\bfOmega$ and $\bfGamma$ are independent.
    Then, $\bfUpsilon : \bbK^{q\times p} \rightarrow \bbK^{k}$ defined in \eqref{equ:def upsilon} is a $(\varepsilon, \delta, d)$ oblivious $\bbK^{q\times p} \rightarrow \bbK^{k}$ subspace embedding.
\end{proposition}
\begin{proof}
    Firstly, \Cref{prop:ose union bound} implies that $\bfX \mapsto \bfX \bfSigma^*$ is a $(\varepsilon / 4, \delta/3, d)$ oblivious $\bbK^{q\times p} \rightarrow \bbK^{q\times k_{\bfSigma}}$ subspace embedding.
    Secondly, \Cref{prop:ose union bound} implies that $\bfY \mapsto \bfY \bfOmega^*$ is a $(\varepsilon / 4, \delta/3, d)$ oblivious $\bbK^{k_{\bfSigma}\times q} \rightarrow \bbK^{k_{\bfSigma}\times k_{\bfOmega}}$ subspace embedding.
    Moreover, since $\|\bfY \bfOmega^*\|_F=\|\bfOmega \bfY^*\|_F$, we have that $\bfY \mapsto \bfOmega \bfY$ is a $(\varepsilon / 4, \delta/3, d)$ oblivious $\bbK^{q \times k_{\bfSigma}} \rightarrow \bbK^{k_{\bfOmega} \times k_{\bfSigma}}$ subspace embedding.
    Finally, \Cref{prop:merge three matrix ose same prop} yields the desired result.
\end{proof}

It is important to note that \Cref{prop:matrix embedding} can be combined with any classical oblivious subspace embeddings $\bfSigma$, $\bfOmega$ and $\bfGamma$.
For example, when the matrix $\bfX$ is only given via vector multiplication from the right, it is interesting to choose $\bfSigma$ as a rescaled Gaussian, and $\bfOmega$ and $\bfGamma$ as structured embeddings which are fast to apply, such as P-SRHT or sparse embeddings.
This allows us to compute $\bfX\bfSigma^*$ through matrix-vector multiplications with $\bfX$, using Gaussian $\bfSigma$ of nearly optimal size, and to efficiently sketch large dense $\bfX\bfSigma^*$ and  $\mathrm{vec}(\bfOmega \bfX\bfSigma^*)$.

A drawback of \Cref{prop:matrix embedding} is that it requires one of the sketching matrix $\bfSigma$ to satisfy the oblivious subspace embedding property with failure probability $\delta/(3q)$, which introduces a dependence on the ambient dimension $q$. 
Although this dependence should have only a minor effect on the sketching dimension $k_{\bfSigma}$, which typically scales just logarithmically with the failure probability, one may still seek guarantees independent of the dimensions of $\bfX$.
The most direct way to obtain such result is to draw $\bfSigma$ from a rescaled Gaussian distribution, and to perform a finer analysis of the concentration of $\|\bfX\bfSigma^*\|_F^2$ than the union bound argument used in \Cref{prop:matrix embedding}.
Indeed, we can show that $\|\bfX \bfSigma^*\|_F^2$ concentrates around $\|\bfX\|_F^2$ at least as fast as $\|\bfSigma \bfx\|_2^2$ concentrates around $\|\bfx\|_2^2$, for $\bfx \in \bbK^p$ and $\bfX\in\bbK^{q\times p}$. 
The first consequence is a refinement of \Cref{prop:matrix embedding} when $\bfSigma$ is drawn from the rescaled Gaussian distribution, as stated in \Cref{prop:matrix embedding gaussian sigma} below.

\begin{proposition}
\label{prop:matrix embedding gaussian sigma}
    Let $\bfSigma \in \bbR^{k_{\bfSigma}\times p}$ be drawn from the rescaled Gaussian distribution.
    Let $\bfOmega \in \bbK^{k_{\bfOmega}\times q}$ be a $(\varepsilon/4, \delta/(3k_{\bfSigma}), d)$ oblivious $\bbK^q \rightarrow \bbK^{k_{\bfOmega}}$ subspace embedding.
    Let $\bfGamma \in \bbK^{k\times k_{\bfSigma} k_{\bfOmega}}$ be a $(\varepsilon/4, \delta/3, d)$ oblivious $\bbK^{k_{\bfSigma} k_{\bfOmega}} \rightarrow \bbK^{k}$ subspace embedding.
    Assume that $\bfSigma$, $\bfOmega$ and $\bfGamma$ are independent.
    If
    \[
        k_{\bfSigma} \geq 168 \varepsilon^{-2} (6.9 \eta_{\bbK} d + \log(6/\delta)),
        \quad \eta_{\bbR} = 1,
        \quad \eta_{\bbC} = 2,
    \]
    then, $\bfUpsilon$ defined in \eqref{equ:def upsilon} is a $(\varepsilon, \delta, d)$ oblivious $\bbK^{q\times p} \rightarrow \bbK^{k}$ subspace embedding.
\end{proposition}
\begin{proof}
    Firstly, \Cref{prop:ose with dim free gaussian} implies that $\bfX \mapsto \bfX \bfSigma^*$ is a $(\varepsilon/4, \delta/3, d)$ oblivious $\bbK^{q\times p} \rightarrow \bbK^{q\times k_{\bfSigma}}$ subspace embedding.
    Then, \Cref{prop:ose union bound} implies that $\bfY \mapsto \bfOmega \bfY$ is a $(\varepsilon / 4, \delta/3, d)$ oblivious $\bbK^{q \times k_{\bfSigma}} \rightarrow \bbK^{k_{\bfOmega} \times k_{\bfSigma}}$ subspace embedding.
    Finally, \Cref{prop:merge three matrix ose same prop} yields the desired result.
\end{proof}

Note that similar results can be derived for $\bfSigma$ with sub-Gaussian random entries, as mentioned in \Cref{rem:hanson wright} below.

\begin{remark}
\label{rem:hanson wright}
    The Gaussian concentration result we used in \Cref{prop:matrix embedding gaussian sigma} is a particular case of the \emph{Hanson-Wright} inequality, see for example \cite[Section 6.2]{vershyninHighDimensionalProbabilityIntroduction2018}, which states a concentration inequality for quadratic forms of independent sub-Gaussian random variables.
    Hence, one may generalize our results to sketching matrices with sub-Gaussian entries, such as Rademacher, to obtain $k_{\bfSigma} = \calO(\varepsilon^{-2} (d + \log(1/\delta)))$, up to modifications in the constants. 
\end{remark}

Note that in \Cref{prop:matrix embedding gaussian sigma},   no dependency of $k_{\bfSigma}$ on the ambient dimensions is introduced.
The second consequence of the Gaussian concentration is that we can derive results for any oblivious embedding $\bfSigma$ without adding dependency to the dimensions of $\bfX$, as depicted in \Cref{prop:matrix embedding gaussian trick}.

\begin{proposition}
\label{prop:matrix embedding gaussian trick}
    Let $\bfSigma\in\bbR^{k_{\bfSigma} \times p}$ (resp. $\bfOmega\in\bbR^{k_{\bfOmega} \times q}$) be a $(\varepsilon', \delta', d)$ oblivious $\bbR^p \rightarrow \bbR^{k_{\bfSigma}}$ (resp. $\bbR^q \rightarrow \bbR^{k_{\bfOmega}}$) subspace embedding, with
    \[
        \varepsilon':= \varepsilon/5,
        \quad 
        \delta' := \delta \lceil 25200 \varepsilon^{-2}(6.9 \eta_{\bbK}d + \log(24/\delta)) \rceil^{-1},
        \quad \eta_{\bbR}=1,
        \quad \eta_{\bbC} = 2.
    \]
    Let $\bfGamma\in\bbR^{k \times k_{\bfSigma} k_{\bfOmega}}$ be a $(\varepsilon/4, \delta/3, d)$ oblivious $\bbR^{k_{\bfSigma} k_{\bfOmega}} \rightarrow \bbR^{k}$ subspace embedding.
    Then, $\bfUpsilon$ defined in \eqref{equ:def upsilon} is a $(\varepsilon, \delta, d)$ oblivious $\bbR^{q\times p} \rightarrow \bbR^k$ subspace embedding.
\end{proposition}
\begin{proof}
    \Cref{prop:ose gaussian trick} implies that $\bfX \mapsto \bfX \bfSigma^*$ (resp. $\bfY \mapsto \bfOmega\bfY$) is a $(\varepsilon/4, \delta/3, d)$ oblivious $\bbR^{q\times p} \rightarrow \bbR^{q\times k}$ (resp. $\bbR^{q\times k_{\bfSigma}} \rightarrow \bbR^{k_{\bfOmega}\times k_{\bfSigma}}$) subspace embedding.
    Then, \Cref{prop:merge three matrix ose same prop} implies the desired result.
\end{proof}

Similarly to \Cref{prop:matrix embedding}, the above \Cref{prop:matrix embedding gaussian trick} can be used with any classical oblivious subspace embeddings $\bfSigma$, $\bfOmega$ and $\bfGamma$.
This can be compared with \Cref{prop:matrix embedding gaussian sigma} which requires $\bfSigma$ to be Gaussian.
The main difference with \Cref{prop:matrix embedding} is that \Cref{prop:matrix embedding gaussian trick} does not introduce additional dependency on the dimension $q$.
In practice, one may use either \Cref{prop:matrix embedding} or \Cref{prop:matrix embedding gaussian trick}, depending on the specific setting.

\subsubsection{Examples of oblivious subspace embeddings}

In this section we provide explicit bounds for the sketch sizes involved in \eqref{equ:def upsilon} for specific choices of random embeddings $\bfSigma$, $\bfOmega$ and $\bfGamma$.
The first case is when drawing these three embeddings from independent rescaled Gaussian distributions, as detailed in \Cref{prop:ose three gaussians}, which yields bounds similar to when using Gaussian embeddings for vectors.
Again, as discussed in \Cref{rem:hanson wright}, one should be able to show that such result holds for general sub-Gaussian distributions, such as the Rademacher distribution, up to modification of the constants.

\begin{proposition}
\label{prop:ose three gaussians}
    Let $\bfSigma \in \bbR^{k \times p}$, $\bfOmega \in \bbR^{k \times p}$ and $\bfGamma \in \bbR^{k \times k^2}$ be independently drawn from rescaled Gaussian distributions, with
    \[
        k \geq 168 \varepsilon^{-2}( 6.9 \eta_{\bbK}d + \log(6/\delta)),
        \quad \eta_{\bbR}=1,
        \quad \eta_{\bbC} = 2.
    \]
    Then, $\bfUpsilon$ defined in \eqref{equ:def upsilon} is a $(\varepsilon, \delta, d)$ oblivious $\bbK^{q\times p} \rightarrow \bbK^k$ subspace embedding.
\end{proposition}
\begin{proof}
    Firstly, \Cref{prop:ose with dim free gaussian} yields that $\bfX \mapsto \bfX \bfSigma^*$ is a $(\varepsilon/4, \delta/3, d)$ oblivious $\bbK^{q\times p} \rightarrow \bbK^{q \times k}$ subspace embedding.
    Similarly, $\bfY \mapsto \bfOmega \bfY$ is a $(\varepsilon/4, \delta/3, d)$ oblivious $\bbK^{q\times k} \rightarrow \bbK^{k \times k}$ subspace embedding, and $\bfGamma$ is a $(\varepsilon/4, \delta/3, d)$ oblivious $\bbK^{k^2} \rightarrow \bbK^k$ subspace embedding.
    Then, \Cref{prop:merge three matrix ose same prop} yields the desired result.
\end{proof}

The problem with Gaussian embeddings is that they can be too costly to use in practice.
Hence, as mentioned below \Cref{prop:matrix embedding}, an efficient way of sketching a matrix $\bfX$, given only by vector multiplication from the right, is to take $\bfOmega$ and $\bfGamma$ as fast to apply embeddings, such as P-SRHT or sparse embeddings.
For example, using \Cref{prop:matrix embedding gaussian sigma} and existing non-asymptotic bounds for the P-SRHT distribution \cite{balabanovRandomizedLinearAlgebra2019}, we obtain the following \Cref{prop:ose gauss srht srht}.

\begin{proposition}
\label{prop:ose gauss srht srht}
    Let $\bfSigma \in \bbR^{k_{\bfSigma} \times p}$ be drawn from the rescaled Gaussian distribution, and let $\bfOmega \in \bbR^{k_{\bfOmega} \times p}$ and $\bfGamma \in \bbR^{k \times k_{\bfSigma} k_{\bfOmega}}$ be drawn from the P-SRHT distribution, with $\bfSigma$, $\bfOmega$ and $\bfGamma$ independent,
    \[
    \begin{aligned}
        k_{\bfSigma} &\geq 168 \varepsilon^{-2}( 6.9 \eta_{\bbK}d + \log(6/\delta)),
        \quad \eta_{\bbR}=1,
        \quad \eta_{\bbC} = 2,
        \\
        k_{\bfOmega} &\geq 32\frac{
            \left(
                \sqrt{d}
                + \sqrt{8 \log(18qk_{\bfSigma} / \delta)}
            \right)^2
        }{
            \varepsilon^2 - \varepsilon^3 / 12
        }
        \log(9dk_{\bfSigma} / \delta),
        \\
        k &\geq 32\frac{
            \left(
                \sqrt{d}
                + \sqrt{8 \log(18k_{\bfOmega}k_{\bfSigma} / \delta)}
            \right)^2
        }{
            \varepsilon^2 - \varepsilon^3 / 12
        }
        \log(9d / \delta).
    \end{aligned}
    \]
    Then, $\bfUpsilon$ defined in \eqref{equ:def upsilon} is a $(\varepsilon, \delta, d)$ oblivious $\bbK^{q\times p} \rightarrow \bbK^{k}$ subspace embedding.
\end{proposition}
\begin{proof}
    Firstly, \Cref{prop:ose with dim free gaussian} implies that $\bfX \mapsto \bfX \bfSigma^*$ is a $(\varepsilon/4, \delta/3, d)$ oblivious $\bbK^{q\times p} \rightarrow \bbK^{q \times k_{\bfSigma}}$ subspace embedding.
    Then, from \cite[Proposition 3.9]{balabanovRandomizedLinearAlgebra2019} the P-SRHT distribution satisfies a $(\varepsilon', \delta',d)$ oblivious $\bbR^q \rightarrow \bbK^{k'}$ subspace embedding property, with $\varepsilon',\delta' \in (0,1)$, if 
    \[
        k' \geq 2\frac{
            \left(\sqrt{d} + \sqrt{\log(8q / \delta')}
            \right)^2
        }{\varepsilon'^2 - \varepsilon'^3/3}
        \log(3d/\delta').
    \]
    Hence, $\bfOmega$ is a $(\varepsilon/4, \delta/(3k_{\bfSigma}), d)$ oblivious $\bbK^q \rightarrow \bbK^{k_{\bfOmega}}$ subspace embedding.
    Similarly, $\bfGamma$ is a $(\varepsilon/4, \delta/3, d)$ oblivious $\bbR^{k_{\bfSigma} k_{\bfOmega}} \rightarrow \bbR^{k}$ subspace embedding.
    Finally, \Cref{prop:matrix embedding gaussian sigma} yields the desired result.
\end{proof}

Let us end this section by discussing on the case $p=q=n \gg 1$.
Then, neglecting logarithmic terms except $\log(1/\delta)$ and $\log\log(n)$, the above \Cref{prop:ose gauss srht srht}, implies that the oblivious embedding property of $\bfUpsilon$ can be ensured with sketch sizes
\[
\begin{gathered}
    k_{\bfSigma} = \calO\left(\varepsilon^{-2}(d + \log(1/\delta))\right),
    \quad
    k_{\bfOmega} = \calO\left(\varepsilon^{-2}(d + \log(1/\delta) + \log(n)) \log(1/\delta)\right),
    \\
    k = \calO\left(\varepsilon^{-2}(d + \log(1/\delta) + \log\log(n))\right).
\end{gathered}
\]
Note that the final sketch size $k$ of $\bfUpsilon$ depends very weakly on the full dimension $n$.
Note also that this dependency can actually be removed using alternative results for the P-SRHT distribution, see \cite{cohenOptimalApproximateMatrix2015}, although the available results are only asymptotic.
Concerning the computational cost, applying P-SRHT embeddings $\bfOmega$ and $\bfGamma$ to a vector costs $\calO(n \log(n))$ flops and $\calO(k_{\bfSigma} k_{\bfOmega} \log(k_{\bfSigma} k_{\bfOmega}))$ flops respectively.
As a result, denoting $T_{\bfX} \geq n\log(n)$ the cost of applying $\bfX$ to a vector, computing $\bfUpsilon(\bfX)\in\bbK^{k_{\bfGamma}}$ costs a number of flops scaling as
\[
    \calO\left(
        \varepsilon^{-2} T_{\bfX} (d + \log(1/\delta))
    \right).
\]

\subsection{Oblivious subspace embeddings in Hilbert-Schmidt spaces}
\label{subsec:ose in hs spaces}

Then, in the same way that we can define oblivious subspace embeddings on $V = (\bbK^n, \innerp{\cdot}{\bfR_V \cdot})$ from oblivious subspace embeddings on $\bbK^n$, we define in \Cref{prop:hs embedding} oblivious subspace embeddings on $HS(V, W)$ from oblivious subspace embeddings on $\bbK^{\dim W \times \dim V}$.

\begin{proposition}
\label{prop:hs embedding}
    Let $V = (\bbK^{\dim V}, \innerp{\cdot}{\bfR_V \cdot})$ and $W = (\bbK^{\dim W}, \innerp{\cdot}{\bfR_W \cdot})$.
    Let $\bfQ_V\in\bbK^{p \times \dim v}$ and $\bfQ_W\in\bbK^{q \times \dim W}$ such that $\bfQ_W^* \bfQ_W = \bfR_W$ and $\bfQ_V^* \bfQ_V = \bfR_V$.
    If $\bfUpsilon$ is a $(\varepsilon, \delta, d)$ oblivious $\bbK^{q\times p} \rightarrow \bbK^k$ subspace embedding, then 
    \begin{equation}
    \label{equ:hs embedding}
        \bfTheta : \bfG \mapsto \bfUpsilon(\bfQ_W \bfG \bfR_V^{-1} \bfQ_V^*)
    \end{equation}
    is a $(\varepsilon, \delta, d)$ oblivious $HS(V, W) \rightarrow \bbK^k$ subspace embedding.
\end{proposition}
\begin{proof}
    This result is obtained with the same reasoning as in \cite[Proposition 3.11]{balabanovRandomizedLinearAlgebra2019}.
    From \eqref{equ:hs as frobenius} we have 
    $\innerp{\cdot}{\cdot}_{HS(V, W)} = \innerp{\bfQ_W \cdot \bfR_V^{-1} \bfQ_V^*}{\bfQ_W \cdot \bfR_V^{-1} \bfQ_V^*}_F$, and from the definition of $\bfTheta$ we have
    $\innerp{\cdot}{\cdot}^{\bfTheta}_{HS(V, W)} = \innerp{\bfQ_W \cdot \bfR_V^{-1} \bfQ_V^*}{\bfQ_W \cdot \bfR_V^{-1} \bfQ_V^*}^{\bfUpsilon}_F$, which yields the desired result.
\end{proof}

Moreover, for low dimensional subspaces $V_p \subset V$ and $W_q \subset W$, we can use the same reasoning to define oblivious embedding on $HS(V_p, W_q)$, as stated in \Cref{prop:hs embedding seminorm}.

\begin{proposition}
\label{prop:hs embedding seminorm}
    Let two subspaces $V_p \subset V = (\bbK^{\dim V}, \innerp{\cdot}{\bfR_V \cdot})$ and $W_q \subset W = (\bbK^{\dim W}, \innerp{\cdot}{\bfR_W \cdot})$.    
    Let $\bfV_p\in\bbK^{\dim V \times p}$ and $\bfW_q\in\bbK^{\dim W \times q}$ whose columns span $V_p$ and $W_q$ respectively with $\bfV_p^* \bfR_V \bfV_p = \bfI_p$ and $\bfW_q^* \bfR_W \bfW_q = \bfI_q$.
    Let $\bfQ_{V_p} = \bfV_p^* \bfR_V$ and $\bfQ_{W_q} = \bfW_q^* \bfR_W$ such that $\bfQ_{V_p}^* \bfQ_{V_p} = \bfR_V \bfPi_{V_p}$ and $\bfQ_{W_q}^* \bfQ_{W_q} = \bfR_W \bfPi_{W_q}$ respectively.
    If $\bfUpsilon$ is a $(\varepsilon, \delta, d)$ oblivious $\bbK^{q\times p} \rightarrow \bbK^k$ subspace embedding, then
    \begin{equation}
    \label{equ:hs embedding seminorm}
        \bfTheta : \bfG \mapsto 
        \bfUpsilon(\bfQ_{W_q} \bfG \bfR_V^{-1} \bfQ_{V_p}^*)
        = \bfUpsilon(\bfW_q^* \bfR_W \bfG \bfV_p)
    \end{equation}
    is a $(\varepsilon, \delta, d)$ oblivious $HS(V_p, W_q) \rightarrow \bbK^k$ subspace embedding.
\end{proposition}
\begin{proof}
    This result is obtained with the same reasoning as in \Cref{prop:hs embedding}.
\end{proof}

We end this section by a few remarks.
Firstly, if in \Cref{prop:hs embedding seminorm} we have $V_p = V$, then we can take $\bfQ_{V_p}$ as in \Cref{prop:hs embedding}.
The same remark applies to the case $W_q = W$.
Secondly, if the dimension $p$ in \Cref{prop:hs embedding seminorm} is rather small, then it is not useful to consider a random embedding $\bfSigma$ in \Cref{prop:matrix embedding}, as one can simply take $\bfSigma = \bfI_p$ without prohibitive cost.
The same remark applies if $q$ is small.
We discuss in more details on the practical aspects that should drive the choice of the sketch in \Cref{sec:practical aspects preconditioners}.

\begin{remark}
\label{rem:sketched mor}
    One may also consider equipping the space $W$ with some sketched inner product $\innerp{\cdot}{\cdot}^{\bfPhi}_W$ for some $W\rightarrow \bbK^k$ embedding $\bfPhi$.
    This would for example allow us leveraging the randomized model reduction framework introduced in \cite{balabanovRandomizedLinearAlgebra2019}.
    This is left to further investigation.
\end{remark}

\begin{remark}[Extension to infinite dimensional setting]
    In the infinite dimensional setting, where both $V$ and $W$ are infinite dimensional Hilbert spaces of functions, a natural extension of the random embedding we proposed would be
    \[
        \bfX \mapsto \bfGamma \mathrm{vec}\big(\innerp{\bfomega_i}{\bfX \bfs_j}_W)_{1\leq i,j\leq k}\big),
    \]
    with $\bfs_1, \cdots, \bfs_k \in V$ and $\bfomega_1, \cdots, \bfomega_k \in W$ random functions drawn from Gaussian distributions with prescribed covariance operators, and with $\bfGamma \in \bbK^{k\times k^2}$.
    There are two main difficulties in this setting.

    First, as pointed out in \cite{perssonRandomizedNystromApproximation2025}, one would want the dominant eigenspaces of the covariance operators to be well aligned, in some sense, with the dominant left or right singular spaces of the operator.
    This alignment impacts the theoretical guarantees of randomized methods extended to the infinite dimensional setting, such as the randomized singular value decomposition in \cite{boulleGeneralizationRandomizedSingular2022,boulleLearningEllipticPartial2023} or the randomized Nystr\"om approximation \cite{perssonRandomizedNystromApproximation2025}.
    This aspect is left to further investigation.

    The second difficulty concerns the measures of quality, since operator norm and HS norm are not equivalent anymore.
    For example, the identity operator is not HS while being bounded in operator norm.
    Hence, it may be critical to consider HS seminorms associated to some finite dimensional reduced space instead, as this allows to get back to the finite dimensional setting.
\end{remark}

\subsection{Comparison with existing methods}
\label{subsec:comparison with existing methods}

In this section we discuss on other random embeddings that have been or could have been used to embed a large matrix using only matrix-vector products.
A first approach could be to consider a structured random sampling embedding applied on the vectorized matrix.
This can be obtained by taking $\bfOmega$ and $\bfSigma$ as random sampling embeddings and $\bfGamma$ as identity.
The resulting embedding $\bfOmega\bfX\bfSigma^*$ then contains the entries $(\bfX)_{i,j}$ for all indices $i$ and $j$ sampled by $\bfOmega$ and $\bfSigma$ respectively.
However, such random embedding is not well suited for oblivious embeddings, as it performs well only when the indices are sampled according to the \emph{leverage scores} of $\bfX$, resulting in a non-oblivious embedding.

Sketching methods for matrices have been used for example for solving large eigenvalue problems.
First in \cite{elsworthConversionsBarycentricRKFUN2019} a simple scalar valued embedding was considered, so that the sketch of $\bfX$ was $\bfomega^* \bfX \bfs \in \bbK$ with random vectors $\bfomega$ and $\bfs$.
A direct extension was proposed in \cite{guttelRandomizedSketchingNonlinear2024}, where a vector valued sketch is constructed as $(\bfomega_i^* \bfX \bfs_i)_{1\leq i\leq k}$ with random vectors $\bfomega_i$ and $\bfs_i$.
A further extension was proposed in \cite{bujanovicSubspaceEmbeddingRandom2025}, by introducing a random embedding $\bfTheta$ defined as the Khatri-Rao product of two random matrices.
The $i$-th column vector of $\bfTheta^*$ is $\mathrm{vec}(\bfomega_i \bfs_i^*)$ for random vectors $\bfomega_i$ and $\bfs_i$, so that 
\[
    \bfTheta \mathrm{vec}(\bfX) = (\bfomega_i^* \bfX \bfs_i)_{1\leq i\leq k} = \mathrm{diag}(\bfOmega \bfX \bfSigma^*).
\]
Importantly, \cite{bujanovicSubspaceEmbeddingRandom2025} provided an oblivious subspace embedding result for this Khatri-Rao embedding, which requires drawing $\bfOmega$ and $\bfSigma$ from the rescaled Gaussian distribution, with
\[ 
    k = \calO (\varepsilon^{-2} d^{3/2} + \varepsilon^{-2} d \log(1/\delta) + d^{1/2} \varepsilon^{-1} \log(1/\delta)^2),
\]
where the constant is much larger than the constant in \Cref{prop:ose three gaussians}.
The first difference compared to our approach is that we provide, in \Cref{prop:matrix embedding} and \Cref{prop:matrix embedding gaussian trick}, oblivious subspace embedding results for arbitrary oblivious subspace embeddings $\bfSigma$ and $\bfOmega$, including efficient embeddings such as P-SRHT or sparse embeddings.
This can be important, for example, when $\bfX$ has a factored form $\bfU \bfV$, so that $\bfOmega \bfX\bfSigma^* = (\bfOmega \bfU)(\bfSigma\bfV)^*$ can be computed efficiently using structured embeddings.

A deeper comparison between this Khatri-Rao embedding result and our results from \Cref{subsec:ose in frobenius space} requires specifying the available operations with $\bfX$.
A first setting, which is the one we are interested in the present work, is when we can compute $\bfX \bfs$ for arbitrary vectors $\bfs$.
Then from \Cref{prop:ose gauss srht srht}, we can compute $\bfUpsilon(\bfX)$ where the dominant computational cost is the $\calO(d)$ matrix-vector product with $\bfX$, fixing $\varepsilon,\delta$.
On the other hand, computing Khatri-Rao embedding requires $\calO(d^{3/2})$ matrix-vector product with $\bfX$. 
Hence, in this first setting, our approach is better suited.

Another setting is when only information on $\bfX$ of the form $\bfomega^* \bfX \bfs$ are available.
Then, computing the Khatri-Rao embedding requires $\calO(d^{3/2})$ operations on $\bfX$, whereas using our approach with \Cref{prop:ose three gaussians} requires $\calO(d^2)$ operations on $\bfX$.
Hence, in this second setting, the Khatri-Rao embedding is better suited.


\section{Measures of quality of a Preconditioner}
\label{sec:measures of quality of a preconditioner}

In this section, we present different measures of quality of a preconditioner in terms of norms described in \Cref{sec:norms on linear operators}.
The choice for the norm will be driven by two main incentives.
The first incentive is the purpose of the preconditioner, as a more specific purpose can leverage weaker norms, which will yield sharper a posteriori estimators.
The second incentive, more practical, is the computability of an approximate solution to \eqref{equ:proj error A onto precond space}.

Firstly in \Cref{subsec:general purpose preconditioner} we present measures of quality of a preconditioner in a general context.
Secondly in \Cref{subsec:mor-purpose preconditioner} we present measures of quality of a preconditioner in a model reduction context.
It is important to keep in mind that all the theoretical results in these first two sections involving the operator norms, which are challenging to minimize, also hold for the Hilbert-Schmidt norms, whose minimizer will be efficiently estimated using random sketching methods.
It is also worth noting that all the results in these first two sections can be applied for any given preconditioner $\bfP$, regardless of the set from which it is taken.
Finally, in \Cref{subsec:sketched measures of quality} we present sketched versions of the measures of quality based on the Hilbert-Schmidt norms, which leverage the fact that $\bfP$, and therefore $\bfE:= \bfI - \bfP\bfA$, lie in a low-dimensional vector space.

\subsection{General purpose preconditioner}
\label{subsec:general purpose preconditioner}

A first choice of norm for evaluating the general quality of a preconditioner is the operator norm defined in \eqref{equ:def operator norm} with $V=W=U$, which leads to the measure of quality $\|\mathbf{E}\|_{U, U}$, with $\bfE = \bfI - \bfP \bfA$.
This is denoted as a \emph{multipurpose} indicator as it can be used to bound the lowest and largest singular values of the preconditioned operator $\bfP \bfA$, assuming that $\|\bfE\|_{U, U} < 1$, as stated in \Cref{prop:multi purpose bound singvals PA} below.

\begin{proposition}
\label{prop:multi purpose bound singvals PA}
    Assuming that $\|\bfE\|_{U, U} < 1$, it holds
    \begin{equation}
    \label{equ:multi purpose bound singvals PA}
        1 - \|\bfE\|_{U, U}
        \leq \alpha(\bfP\bfA)
        \leq \beta(\bfP\bfA)
        \leq 1 + \|\mathbf{E}\|_{U, U},
    \end{equation}
    with $\alpha(\bfP\bfA) := \min_{\bfv \in U \setminus \{0\}} \frac{\|\bfP \bfA \bfv\|_U}{\|\bfv\|_U}$ and $\beta(\bfP\bfA) := \max_{\bfv \in U \setminus \{0\}} \frac{\|\bfP \bfA \bfv\|_U}{\|\bfv\|_U}$.
\end{proposition}
\begin{proof}
    Let $\bfv \in U$.
    By writing $\|\bfP\bfA\bfv\|_U = \|\bfv - (\bfI - \bfP\bfA)\bfv\|_U$ and applying triangle inequalities, we obtain
    \[
        \|\bfv\|_U - \|(\bfI - \bfP\bfA) \bfv\|_U
        \leq \|\bfP\bfA\bfv\|_U 
        \leq \|\bfv\|_U + \|(\bfI - \bfP\bfA) \bfv\|_U.
    \]
    Then by property of the operator norm it holds $\|(\bfI - \bfP\bfA) \bfv\|_U \leq \|\bfE\|_{U,U} \|\bfv\|_U$, which yields the desired result.
\end{proof}

The above \Cref{prop:multi purpose bound singvals PA} has several direct consequences.
A first consequence is that we can bound the condition number of $\bfP\bfA$,
\begin{equation}
\label{equ:multi purpose bound cond PA}
    \mathrm{cond}_2(\bfP\bfA)
    := \frac{\beta(\bfP\bfA)}{\alpha(\bfP\bfA)}
    \leq \frac{1 + \|\bfE\|_{U,U}}{1 - \|\bfE\|_{U,U}}.
\end{equation}
A second consequence is that we can control the quasi-optimality of the preconditioned residual-based error estimator $\|\bfP \bfr(\bfv)\|_U := \|\bfP\bfb - \bfP\bfA \bfv\|_U$ for any $\bfv \in U$, which satisfies
\begin{equation}
\label{equ:multi purpose residual estimator}
    \frac{1}{1 + \|\bfE\|_{U,U}} \|\bfP \bfr (\bfv)\|_U
    \leq \|\bfu - \bfv\|_U
    \leq \frac{1}{1 - \|\bfE\|_{U,U}} \|\bfP \bfr (\bfv)\|_U.
\end{equation}
Finally, it can also be used to bound the quasi-optimality of the preconditioned Galerkin projection as a direct consequence of \Cref{prop:quasi optim precond galerkin}, which we discuss later in \Cref{subsec:mor-purpose preconditioner}.

The main problem with choosing $\|\bfE\|_{U,U}$ as a measure of quality is that its computation involves computing online the largest singular value of a matrix of size $n\times n$, which is prohibitive.
For the same reason, solving \eqref{equ:proj error A onto precond space} with $\|\cdot\| = \|\cdot\|_{U,U}$ is intractable, as it is a large scale non-quadratic, although convex, optimization problem.
This problem can be partially circumvented by considering the Hilbert-Schmidt norm defined in \eqref{equ:def hs norm} with $V=W=U$.
This leads to the second measure of quality $\|\bfE\|_{HS(U,U)}$, which yields the same result as in \Cref{prop:multi purpose bound singvals PA}, as stated in \Cref{coro:multi purpose bound singvals PA hs}.

\begin{corollary}
\label{coro:multi purpose bound singvals PA hs}
    Assuming that $\|\bfE\|_{HS(U, U)} < 1$, it holds
    \begin{equation}
    \label{equ:multi purpose bound singvals PA hs}
        1 - \|\bfE\|_{HS(U, U)}
        \leq \alpha(\bfP\bfA)
        \leq \beta(\bfP\bfA)
        \leq 1 + \|\bfE\|_{HS(U, U)},
    \end{equation}
    with $\alpha(\bfP\bfA)$ and $\beta(\bfP\bfA)$ as defined in \Cref{prop:multi purpose bound singvals PA}.
\end{corollary}
\begin{proof}
    It is a direct consequence of \Cref{prop:multi purpose bound singvals PA} and the right inequality in \Cref{prop:hs norm equivalent op norms}.
\end{proof}

As already mentioned in \Cref{sec:norms on linear operators}, the main advantage of $\|\cdot \|_{HS(U,U)}$ is that we explicitly know the solution to \eqref{equ:proj error A onto precond space}.
Indeed, using \eqref{equ:hs as frobenius} and the fact that $\|\cdot\|_F = \|\mathrm{vec}(\cdot)\|_2$ where $\mathrm{vec}$ is a vectorization operator, solving \eqref{equ:proj error A onto precond space} is equivalent to solving the following linear least-squares problem of size $n^2 \times p$,
\[
    \min_{\bfa \in \bbK^p}
    \|\mathrm{vec}(\bfI) 
    - \sum_{i=1}^p a_i 
    \mathrm{vec}(\bfR_U^{1/2}\bfY_i \bfA\bfR_U^{-1/2})\|_2.
\]

The first problem with $\|\cdot \|_{HS(U,U)}$ is that the matrices $(\bfY_i)_{1\leq i\leq p}$ are in general not sparse, thus the above linear system is dense.
Although the terms of the associated normal equation could be precomputed offline if $\bfA = \bfA(\xi)$ is parameter separable (written as a linear combination of a few known operators), the associated offline cost would be $\calO((np)^2)$, which is highly prohibitive given that $n = \dim U \gg 1$.
This problem is circumvented by the sketching approach described in \Cref{subsec:sketched measures of quality}.

The second problem with $\|\cdot \|_{HS(U,U)}$ is that, in practice, it is often hard to achieve $\|\bfE \|_{HS(U,U)} < 1$, which is required in \Cref{coro:multi purpose bound singvals PA hs}.
This is mainly due to the left inequality in \Cref{prop:hs norm equivalent op norms}, which states that $\|\bfE \|_{U,U}$ and $\|\bfE \|_{HS(U,U)}$ can differ by a factor $n^{1/2}$.
In particular, the Hilbert-Schmidt norm can be highly sensitive to the discretization of the initial operator between infinite dimensional spaces, which might not even be Hilbert-Schmidt. 
Despite this, numerical experiments from \Cref{sec:numerics preconditioners} tends to show that minimizing $\|\bfE \|_{HS(U,U)}$ often yields a good preconditioner even when $\|\bfE \|_{HS(U,U)} \geq 1$.
This second problem can be circumvented in the context of model reduction, as described in \Cref{subsec:mor-purpose preconditioner}.

\subsection{Model reduction purpose preconditioner}
\label{subsec:mor-purpose preconditioner}

In this section, we detail choices of norms on operators that are better suited to the context of model reduction (especially when using Galerkin projections) than the general purpose norms from \Cref{subsec:general purpose preconditioner}.
In this section we consider that we are given an $r$-dimensional subspace $U_r \subset U$, spanned by the columns of an orthogonal matrix $\bfU_r \in \bbK^{n\times r}$.
In \Cref{subsubsec:preconditioned galerkin} we define the preconditioned Galerkin projection on $U_r$ and we show a quasi-optimality result.
In \Cref{subsubsec:residual based error estimator} we introduce and analyze a new residual-based error estimator based on a moderately large dimensional subspace $U_m$, and we compare it with other existing methods.

\subsubsection{Preconditioned Galerkin projection}
\label{subsubsec:preconditioned galerkin}

In this setting, a first important seminorm is $\|\bfE\|_{U_r, U_r}$, as defined in \eqref{equ:def operator hs seminorm}.
Indeed, assuming that $\|\bfE\|_{U_r, U_r} < 1$, the operator $\bfP\bfA$ uniquely defines a Galerkin projection $\bfu_r \in U_r$, which we refer to as the \emph{preconditioned Galerkin projection}, as stated in \Cref{prop:def precond gal}.

\begin{proposition}[Preconditioned Galerkin]
\label{prop:def precond gal}
    Assuming that $\|\bfE\|_{U_r, U_r} < 1$, there exists a unique $\bfu_r \in U_r$ satisfying
    \begin{equation}
    \label{equ:def precond gal}
        \|\bfPi_{U_r}\bfP \bfA (\bfu - \bfu_r)\|_U = 0.
    \end{equation}
\end{proposition}
\begin{proof}
    Since $\bfU_r$ has rank $r$, solving \eqref{equ:def precond gal} is equivalent to solving 
    \[
        \bfU_r^* \bfR_U\bfP \bfA \bfu = \bfU_r^* \bfR_U\bfP\bfA \bfU_r \bfa_r,
    \]
    for $\bfa_r \in \bbK^r$,
    Hence, existence and uniqueness are both ensured if $\bfB_r := \bfU_r^* \bfR_U\bfP\bfA \bfU_r \in \bbK^{r\times r}$ is invertible. 
    Let $\bfa \in \bbK^r$ such that $\bfB_r \bfa = 0$ and let $\bfv = \bfU_r \bfa_r$.
    We then have
    \[
        0 
        = \bfPi_{U_r} \bfP \bfA \bfv
        = \bfPi_{U_r} (\bfI - \bfE) \bfv
        = \bfv  - \bfPi_{U_r} \bfE \bfv 
        = \bfv  - \bfPi_{U_r} \bfE \bfPi_{U_r} \bfv.
    \]
    Using this and the definitions of the operator norm and seminorm yields
    \[
        \|\bfv\|_U = \| \bfPi_{U_r} \bfE \bfPi_{U_r} \bfv\|_U \leq \|\bfPi_{U_r} \bfE \bfPi_{U_r}\|_{U,U} \|\bfv\|_U
        = \|\bfE\|_{U_r,U_r} \|\bfv\|_U.
    \]
    Finally, since we assumed $\|\bfE\|_{U_r,U_r}<1$, the above equation implies $\|\bfv\|_U=0$ thus $\|\bfa\|_2=0$, hence $\bfB_r$ is invertible, which yields the desired result.
\end{proof}

An important point is that verifying whether the assumption $\|\bfE\|_{U_r, U_r} < 1$ holds for a fixed $\bfP$ can be done efficiently.
Indeed, computing $\|\bfE\|_{U_r, U_r}$ only involves computing the largest singular value of the small reduced matrix $\bfB_r = \bfU_r^* \bfR_U \bfP \bfA \bfU_r \in \bbR^{r\times r}$, whose online cost is reasonable.
Moreover, as described in \Cref{sec:practical aspects preconditioners}, if $\bfA = \bfA(\xi)$ is parameter separable for some parameter $\xi$, then so is $\bfB_r(\xi)$ with a reasonable number of affine terms, which are small $r\times r$ matrices.
However, building $\bfP$ by solving \eqref{equ:proj error A onto precond space} remains a challenge.
This can be circumvented by considering instead $\|\bfE\|_{HS(U_r, U_r)}$, which differs from $ \|\bfE\|_{U_r, U_r}$ by at most a factor $r^{1/2}$ in view of \Cref{prop:hs norm equivalent op norms}, and for which solving \eqref{equ:proj error A onto precond space} is equivalent to solving a linear least-squares problem of size $r^2 \times p$.

In the non-preconditioned setting, the quasi-optimality of the Galerkin projection is given by a modified C\'ea's lemma, see for example \cite[Proposition 2.2]{balabanovRandomizedLinearAlgebra2019}.
In the preconditioned setting, a similar result holds, which makes use of $\|\bfE\|_{U_r, U_r}$ as well as another seminorm of interest, $\|\bfE\|_{U, U_r}$.
This is stated in \Cref{prop:quasi optim precond galerkin}.

\begin{proposition}
\label{prop:quasi optim precond galerkin}
    Assume that $\|\bfE\|_{U_r, U_r} < 1$.
    Then for $\bfu_r \in U_r$ the solution to the preconditioned Galerkin system \eqref{equ:def precond gal}, it holds
    \begin{equation}
    \label{equ:quasi optim precond galerkin}
        \|\bfu - \bfu_r\|_U \leq 
        (1 + \frac{\|\bfE\|_{U, U_r}}{1 - \|\bfE\|_{U_r, U_r}})
        \|\bfu - \bfPi_{U_r} \bfu\|_U,
    \end{equation}
    with $\|\cdot\|_{U,U_r}$ and $\|\cdot\|_{U_r, U_r}$ as defined in \eqref{equ:def operator hs seminorm}.
\end{proposition}
\begin{proof}
    By triangle inequality, we first have that
    \[
        \|\bfu - \bfu_r\|_U
        \leq \|\bfu - \bfPi_{U_r}\bfu\|_U + \|\bfPi_{U_r} \bfu - \bfu_r\|_U.
    \]
    Let us now adequately bound $\|\bfPi_{U_r} \bfu - \bfu_r\|_U$, which is equal to $\|\bfPi_{U_r} (\bfu - \bfu_r)\|_U$ as $\bfu_r \in U_r$.
    Noting that $\bfI = \bfE + \bfP\bfA$ and using the triangle inequality and the fact that $\bfu_r$ satisfies \eqref{equ:def precond gal}, we obtain
    \[
        \|\bfPi_{U_r} (\bfu - \bfu_r)\|_U
        \leq \|\bfPi_{U_r} \bfE (\bfu - \bfu_r)\|_U
        + \|\bfPi_{U_r} \bfP\bfA (\bfu - \bfu_r)\|_U
        =  \|\bfPi_{U_r} \bfE (\bfu - \bfu_r)\|_U.
    \]
    Then, using the triangle inequality, the fact that $\bfu_r \in U_r$, the definition of the operator norms from \eqref{equ:def operator norm} and the definition of the associated seminorms from \eqref{equ:def operator hs seminorm}, we obtain
    \[  
    \begin{aligned}
        \|\bfPi_{U_r} (\bfu - \bfu_r)\|_U
        & \leq \|\bfPi_{U_r} \bfE (\bfu - \bfPi_{U_r}\bfu)\|_U 
        + \|\bfPi_{U_r} \bfE (\bfPi_{U_r}\bfu - \bfu_r)\|_U \\
        & = \|\bfPi_{U_r} \bfE (\bfu - \bfPi_{U_r}\bfu)\|_U 
        + \|\bfPi_{U_r} \bfE \bfPi_{U_r} (\bfu - \bfu_r)\|_U \\
        & \leq \|\bfPi_{U_r}\bfE\|_{U, U} \|\bfu - \bfPi_{U_r}\bfu\|_U + \|\bfPi_{U_r}\bfE\bfPi_{U_r}\|_{U, U} \|\bfu - \bfu_r\|_U \\
        & = \|\bfE\|_{U, U_r} \|\bfu - \bfPi_{U_r}\bfu\|_U + \|\bfE\|_{U_r, U_r} \|\bfu - \bfu_r\|_U.
    \end{aligned}
    \]
    Finally, combining this last inequality with the first inequality of the proof, the assumption $\|\bfE\|_{U_r, U_r}<1$ yields the desired inequality.
\end{proof}

The problem with the above \Cref{prop:quasi optim precond galerkin} is that computing $\|\bfE\|_{U, U_r}$ for a fixed $\bfP$ is relatively challenging, as it requires computing the largest singular value of the wide matrix $\bfU_m^* \bfR_U \bfP \bfA \in \bbK^{r \times n}$.
Moreover, building $\bfP$ by solving \eqref{equ:proj error A onto precond space} is even more challenging.
The first step to circumvent this is to consider instead $\|\bfE\|_{HS(U,U_r)}$, which differs from $\|\bfE\|_{U,U_r}$ by at most a factor $r^{1/2}$ in view of \Cref{prop:hs norm equivalent op norms}, and for which solving \eqref{equ:proj error A onto precond space} is equivalent to solving a linear least-squares problem of size $nr\times p$.
The problem is that the associated cost depends on the full dimension $n$, which is not online efficient.
One way to circumvent this could be to use the associated normal equation, which inherits parameter separability from $\bfA = \bfA(\xi)$, and which would also allow to efficiently compute $\|\bfE\|_{U, U_r}$ online.
The problem is that the corresponding number of affine terms is essentially squared, which may result in prohibitive offline cost as well as numerical instability due to round-off errors.
This problem can be circumvented by using instead a sketched norm $\|\bfE\|_{HS(U, U_r)}^{\bfTheta}$ as defined in \Cref{sec:random sketching for linear operators,subsec:sketched measures of quality} and detailed in \Cref{sec:practical aspects preconditioners}.

\begin{remark}
\label{rem:minimizing weighted sum}
    In addition to computing the quasi-optimality constant $1 + \frac{\|\bfE\|_{U, U_r}}{1 - \|\bfE\|_{U_r, U_r}}$ from \eqref{equ:quasi optim precond galerkin}, one may want to minimize it.
    A first approach is to use the fact that $\|\bfE\|_{U_r, U_r} \leq \|\bfE\|_{U, U_r}$ and only focus on minimizing $\|\bfE\|_{U, U_r}$.
    A second approach is to minimize instead a weighted sum $\|\bfE\|_{U, U_r}^2 + \tau^2 \|\bfE\|_{U, U_r}^2$ for some $\tau >0$.
    In particular, calculations show that 
    \[
        \|\bfE\|_{U, U_r}^2 
        + \tau^2 \|\bfE\|_{U, U_r}^2 
        \leq \frac{1}{2}
        \implies 
        \frac{\|\bfE\|_{U, U_r}}{1 - \|\bfE\|_{U_r, U_r}}
        \leq \sqrt{2} \tau.
    \]
    As a result, the user can define some prescribed bound on the quasi-optimality constant, define $\tau$ accordingly and then minimize $\|\bfE\|_{U, U_r}^2 + \tau^2 \|\bfE\|_{U, U_r}^2$.
    In practice, minimizing instead $\|\bfE\|_{HS(U, U_r)}^2 + \tau^2 \|\bfE\|_{HS(U, U_r)}^2$ is equivalent to solving a linear least-squares problem.
\end{remark}

\subsubsection{Residual-based error estimator}
\label{subsubsec:residual based error estimator}

Now when it comes to a posteriori error estimators based on the preconditioned residual $\|\bfP\bfr(\bfu_r)\|_U = \|\bfP\bfb - \bfP \bfA \bfu_r\|_U$, the available results are much less satisfying.
Indeed, the seminorms $\|\bfE\|_{U, U_r}$ and $\|\bfE\|_{U_r, U_r}$ are not sufficient anymore in order to ensure a result similar to \eqref{equ:multi purpose residual estimator}.
What we can do instead is to project the preconditioned residual onto some test space $U_m$ with moderately large $m$, such that $r\leq m\ll n$ and $U_r \subset U_m$.
We will then use $\|\bfPi_{U_m} \bfP \bfr(\bfu_r)\|_U$ as an a posteriori error indicator in order to estimate $\|\bfu - \bfu_r\|_U$.
Under the right assumptions, \Cref{prop:bound error by precond residual} states that this estimator can be used to both lower and upper bound the approximation error using the seminorms $\|\bfE\|_{U, U_m}$ and $\|\bfE\|_{U_m, U_m}$.
Note that since $U_r \subset U_m$, we have that $\|\bfE\|_{U_r, U_r} \leq \|\bfE\|_{U_m, U_m}$ and $\|\bfE\|_{U, U_r} \leq \|\bfE\|_{U, U_m}$, thus it is also relevant regarding the quasi-optimality of the preconditioned Galerkin projection to minimize $\|\bfE\|_{U_m, U_m}$, $\|\bfE\|_{U, U_m}$ or a weighted sum as in \Cref{rem:minimizing weighted sum}.

\begin{proposition}
\label{prop:bound error by precond residual}
    Let $\bfv \in U_m$.
    Assume that there exists $\tau \in [0, 1)$ such that 
    \[
    \|\bfu - \bfPi_{U_m} \bfu \|_U \leq 
    \tau \|\bfu - \bfv\|_U
    \quad \text{and} \quad
    \|\bfE\|_{U_m, U_m} + \tau \|\bfE\|_{U, U_m} < \sqrt{1-\tau^2}.
    \]
    Then,
    \begin{equation}
    \label{equ:bound error by precond residual}
        \frac{
            \|\bfPi_{U_m} \bfP \bfr(\bfv)\|_U
            }{
            1 + \|\bfE\|_{U_m, U_m} + \tau \|\bfE\|_{U, U_m}
            }
        \leq \|\bfu - \bfv\|_U
        \leq 
        \frac{
            \|\bfPi_{U_m} \bfP \bfr(\bfv)\|_U
            }{
            \sqrt{1-\tau^2} - \|\bfE\|_{U_m, U_m} - \tau \|\bfE\|_{U, U_m}
            }
    \end{equation}
\end{proposition}
\begin{proof}
    First, by projection property and the fact that $\bfv \in U_m$, we obtain
    \[
        \|\bfPi_{U_m} (\bfu - \bfv)\|_U^2
        \leq \|\bfu - \bfv\|_U^2
        = \|\bfPi_{U_m}(\bfu - \bfv)\|_U^2 + \|\bfu - \bfPi_{U_m} \bfu\|_U^2.
    \]
    Then, by using the first assumption we obtain
    \[
        \|\bfu - \bfPi_{U_m} \bfu\|_U^2
        \leq \tau^2 \|\bfu - \bfv\|_U^2
        \leq \frac{\tau^2}{1 - \tau^2} \|\bfPi_{U_m}(\bfu - \bfv)\|_U^2.
    \]
    Gathering the last two equations yields
    \begin{equation}
    \label{equ:proof precond residual 1}
        \|\bfPi_{U_m} (\bfu - \bfv)\|_U^2
        \leq \|\bfu - \bfv\|_U^2
        \leq \frac{1}{1 - \tau^2} 
        \|\bfPi_{U_m} (\bfu - \bfv)\|_U^2.
    \end{equation}
    Now, let us prove the right inequality in \eqref{equ:bound error by precond residual}.
    Using two triangle inequalities, the fact that $\bfv \in U_m$, the definition of the operators seminorms and the first assumption on $\tau$, we obtain 
    \[
    \begin{aligned}
        \|\bfPi_{U_m} (\bfu - \bfv)\|_U
        & \leq \|\bfPi_{U_m} \bfP \bfr(\bfv)\|_U
        + \|\bfPi_{U_m} \bfE (\bfu - \bfv)\|_U \\
        & \leq \|\bfPi_{U_m} \bfP \bfr(\bfv)\|_U
        + \|\bfPi_{U_m} \bfE \bfPi_{U_m} (\bfu - \bfv)\|_U 
        + \|\bfPi_{U_m} \bfE (\bfu - \bfPi_{U_m} \bfu)\|_U  \\
        & \leq \|\bfPi_{U_m} \bfP \bfr(\bfv)\|_U
        + \|\bfE\|_{U_m, U_m} \|(\bfu - \bfv)\|_U 
        + \|\bfE \|_{U, U_m} \|\bfu - \bfPi_{U_m} \bfu\|_U  \\
        &  \leq \|\bfPi_{U_m} \bfP \bfr(\bfv)\|_U
        + (\|\bfE\|_{U_m, U_m} + \tau \|\bfE\|_{U, U_m} )\|(\bfu - \bfv)\|_U.
    \end{aligned}
    \]
    Then from \eqref{equ:proof precond residual 1} we have $\sqrt{1-\tau^2} \|\bfu - \bfv\|_U \leq \|\bfPi_{U_m} (\bfu - \bfv)\|_U$, which combined with the previous inequality yields the desired right inequality in \eqref{equ:bound error by precond residual}.
    For the left inequality, starting by a reverse triangle inequality and following the same reasoning, we obtain
    \[
        \|\bfPi_{U_m}(\bfu - \bfv)\|_U 
        \geq \|\bfPi_{U_m} \bfP \bfr(\bfv)\|_U
        -  (\|\bfE\|_{U_m, U_m} + \tau \|\bfE\|_{U, U_m} )\|(\bfu - \bfv)\|_U.
    \]
    Then from \eqref{equ:proof precond residual 1} we have $\|\bfPi_{U_m} (\bfu - \bfv)\|_U \leq\|\bfu - \bfv\|_U$, which combined with the previous inequality yields the desired left inequality in \eqref{equ:bound error by precond residual}.
\end{proof}

The main advantage of our error estimator compared to classical approaches, such as the successive constraint method \cite{huynhSuccessiveConstraintLinear2007}, is that we do not approximate the inf-sup constant of a high-dimensional parameter-dependent operator, whose online cost can be expensive and which can yield pessimistic error bounds.
Indeed, our estimator can be efficiently estimated online, as long as we have access to some $\bfU_m\in\bbK^{n \times m}$ whose columns are an orthonormal basis of $U_m$.
Then, $\|\bfPi_{U_m}\bfP\bfr(\bfu_r)\|_U = \|\bfU_m^* \bfR_U \bfP\bfr(\bfu_r)\|_U$ and $\bfU_m^* \bfR_U \bfP\bfr(\bfu_r)$ directly inherit from parameter separability of $\bfA$ and $\bfb$.
Other related works, especially \cite{hainHierarchicalPosterioriError2019,smetanaRandomizedResidualBasedError2019}, share this first advantage.

In \cite{hainHierarchicalPosterioriError2019} a hierarchical error estimator was proposed, defined as the distance between the Galerkin projections on $U_m$ and $U_r$.
The main advantage of our approach compared to this one is that we do not compute online approximation in $U_m$, which can be very costly if $m$ is too large.
In particular, as described in \Cref{sec:practical aspects preconditioners}, the online costs for estimating $\|\bfE\|_{HS(U_m, U_m)}$ and $\|\bfE\|_{HS(U, U_m)}$ depend very weakly on $m$.
We can also use an approach similar to \cite{balabanovRandomizedLinearAlgebra2019} to efficiently approximate our residual-based estimator on finite parameter sets with offline cost independent on $m$.

Interestingly, the analysis of the estimator in \cite{hainHierarchicalPosterioriError2019} relies on a \emph{saturation} assumption, typically found in hierarchical methods, see for example \cite{bankPosterioriErrorEstimates1993,wohlmuthHierarchicalPosterioriError1999,huangNewPosterioriError2011}.
They assumed that $\|\bfu - \bfu_m \|_U \leq \tau \|\bfu - \bfu_r\|_U$ for some $\tau \in (0,1)$, which is actually slightly more restrictive than our assumption $\|\bfu - \bfPi_{U_m} \bfu \|_U \leq \tau \|\bfu - \bfu_r\|_U$ in \Cref{prop:bound error by precond residual}.
Note that in \cite{hainHierarchicalPosterioriError2019} the authors investigated more deeply this assumption as well as the construction of $U_m$, and proposed an offline estimation of $\tau$.
This could be beneficial for our approach, as we do not provide yet any method for verifying this assumption a posteriori. 
This is left to further investigation.
We mention that numerical experiments from \Cref{sec:numerics preconditioners} tends to show that our error estimator is accurate 

In \cite{smetanaRandomizedResidualBasedError2019} a randomized error estimator was proposed, defined by projecting the residual onto approximate solutions to dual problems $\bfA(\xi)^* \bfy_i(\bfxi)=\bfz_i$ with random right-hand sides $\bfz_i$ for $1\leq i\leq k$.
These dual solutions are approximated using a reduced basis approach computed with Galerkin projection onto a space $Y_m$.
Their approach can be seen as a randomized approximation of the preconditioned residual norm $\|\bfP \bfr(\bfv)\|_{U}$ from \eqref{equ:multi purpose residual estimator}, with $\bfP(\xi) = \bfY_m (\bfY_m^* \bfA(\xi) \bfY_m)^{-1} \bfY_m^*$ where $\bfY_m$ is a matrix whose columns form an orthonormal basis of $Y_m$.

There are two main advantages of our approach compared to the one in \cite{smetanaRandomizedResidualBasedError2019}.
The first is that we can expect the approximation of $\bfA(\bfxi)^{-1}$ to be better with respect to seminorms tailored to a reduced space, as we proposed, than with respect to general norms, as in \cite{smetanaRandomizedResidualBasedError2019} and \Cref{subsec:general purpose preconditioner}.
The second is that our linear approximation of $\bfA(\xi)^{-1}$ allows us to obtain guarantees of our randomized methods with high probability on non-finite parameter sets if $\bfA(\bfxi)$ is parameter separable, as detailed in \Cref{subsec:sketched measures of quality}.
This is to be compared to \cite{smetanaRandomizedResidualBasedError2019} in which guarantees only hold for finite parameter sets.
However, it is important to note that the analysis in \cite{smetanaRandomizedResidualBasedError2019} does not rely on the aforementioned \emph{saturation} assumption.

\subsection{Sketched measures of quality}
\label{subsec:sketched measures of quality}

In this section, we describe the randomized (or sketched) measures of qualities that we will use to efficiently construct preconditioners.
In particular, we leverage the fact that the error matrix lies in a low dimensional parameter dependent subspace of linear operators,
\[
    \bfE
    \in \spanv{\bfI, \bfY_1\bfA, \cdots, \bfY_p \bfA},
\]
whose dimension is at most $p+1$.
As a result, for fixed $V\in \{U, U_m\}$, $W \in \{U, U_m\}$, if $\bfTheta$ is a $(\varepsilon, \delta, l)$ oblivious $HS(V, W) \rightarrow \bbK^k$ subspace embedding as defined in \Cref{def:data oblivious subspace embedding}, then it holds
\begin{equation}
\label{equ:embedding property sketched measures}
    \forall \bfP\in\spanv{\bfY_1, \cdots, \bfY_p}, ~
    \big|
    \|\bfE\|_{HS(V,W)}^2 
    - (\|\bfE\|_{HS(V,W)}^{\bfTheta})^2
    \big|
    \leq \varepsilon \|\bfE\|_{HS(V,W)}^2
\end{equation}
with probability at least $1-\delta$ and $l = p+1$.
Recall that, from \Cref{prop:ose gauss srht srht}, we can efficiently sketch a high-dimensional operator at cost of applying this operator to $\calO(\varepsilon^{-2}(l + \log(1/\delta)))$ Gaussian vectors.
Let us now discuss on extensions of the embedding property \eqref{equ:embedding property sketched measures} in the parameter dependent setting, using the same approach as in \cite{balabanovRandomizedLinearAlgebra2019,balabanovRandomizedLinearAlgebra2021}.

Firstly, we may want \eqref{equ:embedding property sketched measures} to hold simultaneously for all operators $\bfA(\xi)$, $\bfxi \in\calP$, with high probability.
There are two ways to ensure that.
Firstly if $\#\calP < +\infty$, then from a union bound of probability argument, taking $\delta = \delta^* / \#\calP$ implies that \eqref{equ:embedding property sketched measures} holds for all $\xi\in\calP$ with probability at least $1-\delta^*$.
Secondly if $\bfA$ is parameter separable such that $\bfA(\xi) = \sum_{j=1}^{m_{\bfA}} \theta^{\bfA}_j(\xi) \bfA_j$, then the error matrix lies in a low dimensional parameter independent subspace of linear operators,
\[
    \forall \xi \in \calP, \quad 
    \bfE(\xi) 
    \in \spanv{\bfI, \bfY_1\bfA_1, \cdots, \bfY_1\bfA_{m_{\bfA}}, \cdots, \bfY_p \bfA_{m_{\bfA}}}
\]
whose dimension is at most $pm_{\bfA}+1$.
Hence, taking $l = p m_{\bfA}+1$ yields \eqref{equ:embedding property sketched measures} for all $\xi\in\calP$ with probability at least $1-\delta$.

Secondly, we may want \eqref{equ:embedding property sketched measures} to hold for any $p$-dimensional subspace of preconditioners $\spanv{\bfY_1, \cdots, \bfY_p}$. 
This is especially the case when using adaptive algorithms, such as the greedy algorithm considered in \Cref{sec:practical aspects preconditioners}, leveraging the sketched quantities detailed in the current section.
One way to ensure that is to take $(\bfY_i)_{1\leq i\leq p}$ from a finite set of size $M$, then from a union bound of probability argument, taking $\delta = \delta^* / \binom{M}{p}$ implies that \eqref{equ:embedding property sketched measures} holds for all $(\bfY_i)_{1\leq i\leq p}$ with probability at least $1-\delta^*$.
Additionally, if $M = \#\calP$ then taking $\delta = \delta^* / (\binom{M}{p} M )$ implies that \eqref{equ:embedding property sketched measures} holds for all $(\bfY_i)_{1\leq i\leq p}$ and all $\bfxi \in \calP$ with probability at least $1-\delta^*$.
It is important to note that the size of the random matrices scales logarithmically with $\delta$, thus the effect of the binomial factor remains reasonable as $\log (1/\delta) \sim p \log(M/\delta^*)$.

We end this section by emphasizing the three embeddings on HS spaces that we will leverage in \Cref{sec:practical aspects preconditioners}.
We assume that we have access to $\bfQ_U \in \bbK^{s\times n}$ such that $\bfQ_U^* \bfQ_U = \bfR_U$ and to $\bfU_m \in \bbK^{n\times m}$ whose columns span $U_m$ and such that $\bfU_m^* \bfR_U \bfU_m = \bfI$.
Using  \Cref{prop:hs embedding,prop:hs embedding seminorm}, we introduce embeddings $\bfPsi$, $\bfXi$ and $\bfLambda$, which are $HS(U,U) \rightarrow \bbK^{k_{\bfPsi}}$, $HS(U,U_m) \rightarrow \bbK^{k_{\bfXi}}$ and $HS(U_m,U_m) \rightarrow \bbK^{k_{\bfLambda}}$  embeddings respectively, defined by
\begin{equation}
\label{equ:def embedding used in practice}
\begin{aligned}
    \bfPsi (\bfE) &:=  
    \bfUpsilon_{s,s}(\bfQ_U \bfE \bfR_U^{-1} \bfQ^*_U), \\
    \bfXi (\bfE) &:= 
    \bfUpsilon_{m,s} (\bfU_m^*\bfR_U\bfE\bfR_U^{-1}\bfQ^*_U), \\
    \bfLambda (\bfE) &:= 
    \bfUpsilon_{m,m}(\bfU_m^*\bfR_U \bfE \bfU_m), \\
\end{aligned}
\end{equation}
where $\bfUpsilon_{s,s}$, $\bfUpsilon_{m,s}$ and $\bfUpsilon_{m,m}$ are embeddings from \Cref{prop:matrix embedding} on $\bbK^{s\times s}$, $\bbK^{m\times s}$ and $\bbK^{m\times m}$ respectively.
We will then approximate $\|\bfE\|_{HS(U, U)}$, $\|\bfE\|_{HS(U, U_m)}$ and $\|\bfE\|_{HS(U_m, U_m)}$ by respectively
$\|\bfE\|_{HS(U, U)}^{\bfPsi} = \|\bfPsi(\bfE)\|_2$, $\|\bfE\|_{HS(U, U_m)}^{\bfXi} = \|\bfXi(\bfE)\|_2$ and $\|\bfE\|_{HS(U_m, U_m)}^{\bfLambda} = \|\bfLambda(\bfE)\|_2$.
As we will see in the next \Cref{sec:practical aspects preconditioners}, these three embedded quantities will inherit parameter separability from $\bfE$, which will allow efficient offline-online decomposition.


\section{Practical aspects}
\label{sec:practical aspects preconditioners}

In this section, we discuss on the practical aspects for manipulating the various quantities involved in \Cref{sec:measures of quality of a preconditioner} when \eqref{equ:def fom} is parameter dependent.
In particular, we describe how to obtain an efficient offline-online decomposition with reasonable offline and online costs and robustness to numerical instabilities.
We also discuss on how to choose the random matrices used in \Cref{prop:matrix embedding}.
In this section, we assume that the operator $\bfA(\xi)$ and the right-hand side $\bfb(\xi)$ are parameter separable, such that 
\begin{equation}
    \bfA(\xi) = \sum_{j=1}^{m_{\bfA}} 
    \theta_j^{\bfA}(\xi) \bfA_j,
    \quad
    \bfb(\xi) = \sum_{j=1}^{m_{\bfb}} 
    \theta_j^{\bfb}(\xi) \bfb_j.
\end{equation}
We recall that parameter separability can be either naturally deduced from the initial problem itself, or obtained from approximation methods such as the empirical interpolation method \cite{madayGeneralMultipurposeInterpolation2009}.

In the following sections, we will analyze the computational costs for computing various quantities in terms of flops.
The main costs will come from high-dimensional matrix-vector multiplications, which depends on the type of architecture, algorithms and precision used in practice.
As a result, one may want to perform a complexity analysis depending on the specific context.
Still, we can make some general cost analysis under some rather classical assumptions.
More precisely, we will denote $T_{\bfY}$, typically $T_{\bfY}\geq n\log(n)$, the cost, in flops, to apply $\bfY_i$ or its adjoint to a vector, and we will assume that applying to a vector $\bfA_j$, $\bfR_U$, $\bfR_U^{-1}$ and $\bfQ_U$, as well as the corresponding adjoint matrices, takes $\calO(n \log (n))$ flops.
This is justified by the fact that $\bfA_j$ and $\bfR_U$ are parameter independent sparse matrices.

\subsection{Sketched measures of quality}
\label{subsec:sketched measures of quality practical}

For $\|\cdot\|^{\bfTheta}_{HS} \in \{\|\cdot\|_{HS(U,U)}^{\bfPsi}, \|\cdot\|_{HS(U, U_m)}^{\bfXi}, \|\cdot\|_{HS(U_m, U_m)}^{\bfLambda}\}$ and for $\bfP = \sum_{j=1}^{P} a_j \bfY_j$, the sketched measures of quality of $\bfP$ based on HS norms can be written as
\[
    \|\bfI - \bfP \bfA(\xi)\|^{\bfTheta}_{HS}
    =\|\bfh - \bfW(\xi) \bfa \|_2,
    \quad
    \bfh = \bfTheta(\bfI),
    \quad
    \bfW^{(i)}(\xi) = \sum_{i=1}^{m_{\bfA}} 
    \theta_i^{m_{\bfA}}(\xi) \bfTheta(\bfY_j \bfA_i),
\]
where $\bfW^{(i)}(\xi) \in \bbK^k$ denotes the $i$-th column of $\bfW(\bfxi)$.
In other words, the columns of the small matrix $\bfW(\xi) \in \bbK^{k\times p}$ are parameter separable, with as many affine terms as $\bfA(\bfxi)$.

\subsubsection{Online costs}
From the above statements, the online cost for evaluating $\bfW(\xi)$, and therefore $\|\bfI - \bfP\bfA(\xi)\|^{\bfTheta}_{HS}$, is $\calO(kpm_{\bfA})$ flops.
An important consequence is that solving \eqref{equ:proj error A onto precond space} is equivalent to solving the small parameter separable linear least-squares problem of size $k\times p$,
\[
    \min_{\bfa \in \bbK^p}
    \|\bfh - \bfW(\xi) \bfa \|_2,
\]
which can be done with stable methods such as SVD or QR, resulting in an online cost of $\calO(kp^2)$ flops.
Note that for computation (not minimization) of $\|\bfI - \bfP\bfA(\xi)\|^{\bfTheta}_{HS}$ over some finite set $\calP$, we can use a similar approach as \cite{balabanovRandomizedLinearAlgebra2021} by introducing a small $(\varepsilon, \delta / \#\calP, 1)$ oblivious $\bbK^k \rightarrow \bbK^{k'}$ subspace embedding with $k' \leq k$ rows, so that given the solution $\bfa(\bfxi)\in\bbK^p$ to the previous least-squares problem, the online cost for computing $\|\bfPhi\bfh - \bfPhi \bfW(\xi) \bfa(\bfxi)\|_2$ is only $\calO(k'p m_{\bfA})$ flops.

\subsubsection{Offline costs}

The offline cost mainly depends on the cost of evaluating $\bfTheta(\bfY_i \bfA_j)$ for $1\leq i\leq p$ and $1\leq j\leq m_{\bfA}$.
There are essentially two main approaches for evaluating the latter, each one being more efficient in specific contexts.

The first approach consists in taking $\bfGamma$ as a fast to apply matrix such as P-SRHT or a sparse embedding, taking $\bfOmega$ and $\bfSigma$ as Gaussian matrices and observing that for some $\bfQ_V$ and $\bfQ_W$ as in \Cref{prop:hs embedding seminorm} we can write
\[
    \bfTheta(\bfY_i \bfA_j)
    = \bfGamma \mathrm{vec}\big(
        (\bfY_i^*\bfQ_V^*\bfOmega^*)^* 
        (\bfA_j \bfR_U^{-1} \bfQ_W^* \bfSigma^*)
        \big).
\]
In the above equation, computing $\bfY_i^*\bfQ_V^*\bfOmega^* \in \bbK^{n \times k}$ for $1 \leq i\leq p$ at a cost of $\calO(T_{\bfY} k p)$ flops can be done independently of computing $\bfA_j \bfR_U^{-1} \bfQ_W^* \bfSigma^* \in \bbK^{n \times k}$ for $1\leq j\leq m_{\bfA}$ at a cost of $\calO(n \log(n) k m_{\bfA})$ flops.
Then, computing the products of those matrices for all $1\leq i\leq p$ and $1\leq j\leq m_{\bfA}$ costs $\calO(nk^2 p m_{\bfA})$ flops.
Finally applying $\bfGamma$ to $p m_{\bfA}$ vectors of dimension $k^2$ costs $\calO(k^2 \log(k) p m_{\bfA})$, which is negligible compared to the previous costs.
Assuming that $kp\gtrsim \log(n)$, the resulting total offline cost is then 
\[
    \calO(k p (T_{\bfY} + nkm_{\bfA}))
\]
flops.
This approach is well suited to highly parallel architectures as computing $\bfY_i^*\bfQ_V^*\bfOmega^*$ and $\bfA_j \bfR_U^{-1} \bfQ_W^* \bfSigma^*$ can be done in parallel using seeded random numbers generator for storing the random matrices implicitly.
On the other hand it requires maintaining in memory $(p+m_{\bfA})$ large matrices of size $n\times k$, which may be prohibitive if the memory cost is not negligible.
Note that one may additionally consider a similar approach as in \cite{balabanovRandomizedLinearAlgebra2019} and approximate $(\bfY_i^*\bfQ_V^*\bfOmega^*)^* (\bfA_j \bfR_U^{-1} \bfQ_W^* \bfSigma^*)$ by $(\bfPhi \bfY_i^*\bfQ_V^*\bfOmega^*)^* (\bfPhi \bfA_j \bfR_U^{-1} \bfQ_W^* \bfSigma^*)$ where $\bfPhi$ is some fast to apply sketching matrix, such as P-SRHT or a sparse embedding, which would mitigate the $\calO(nk^2)$ cost.
We leave this to further investigation.

The second approach consists in taking $\bfGamma$ and $\bfSigma$ as fast to apply matrices, taking $\bfOmega$ as a Gaussian matrix or identity if $\bfQ_W$ has few rows, and observing that we can write
\[
    \bfTheta(\bfY_i \bfA_j)
    = \bfGamma \mathrm{vec}\big(
        (\bfSigma \bfQ_V \bfR_U^{-1} \bfA_j^*
        \bfY_i^*\bfQ_W^*\bfOmega^*)^*
        \big).
\]
In the above equation, for fixed $1\leq i\leq p$ and $1\leq j\leq m_{\bfA}$, we can then evaluate right to left $\bfSigma \bfQ_V \bfR_U^{-1} \bfA_j^* \bfY_i^*\bfQ_W^*\bfOmega^*$, at a cost depending on $V$ and $W$, as we discuss in the next paragraph.
This approach is well suited to slightly parallel (or nonparallel) architectures with memory constraints, as it does not require storing more than $1$ to $k$ high dimensional vectors.
Note that one can still leverage parallelization over the columns of $\bfOmega^*$.
This second approach will be the one we focus on in the rest of this section as well as the one used in \Cref{sec:numerics preconditioners}.

\paragraph{Sketched $HS(U,U)$ norm.}
We take $\bfGamma$ and $\bfSigma$ as fast to apply matrices, and $\bfOmega$ as a Gaussian or Rademacher matrix.
We then evaluate the terms right to left in 
\[
    \bfPsi(\bfY_i \bfA_j)
    = \bfGamma \mathrm{vec}\big(
        (\bfSigma \bfQ_U \bfR_U^{-1} \bfA_j^*
        \bfY_i^*\bfQ_U^*\bfOmega^*)^*
        \big).
\]
As a result, the total offline cost for $\|\bfE\|_{HS(U,U)}^{\bfPsi}$ is $\calO(kp(T_{\bfY} + m_{\bfA} n\log(n)))$ flops.

\paragraph{Sketched $HS(U,U_m)$ norm.}
We take $\bfGamma$ and $\bfSigma$ as fast to apply matrices, and $\bfOmega$ as a Gaussian, Rademacher or identity matrix.
We then evaluate the terms right to left in 
\[
    \bfXi(\bfY_i \bfA_j)
    = \bfGamma \mathrm{vec}\big(
        (\bfSigma \bfQ_U \bfR_U^{-1} \bfA_j^*
        \bfY_i^*\bfU_m\bfOmega^*)^*
        \big).
\]
As a result, the total offline cost for $\|\bfE\|_{HS(U,U)}^{\bfPsi}$ is $\calO(p (k \wedge m) (T_{\bfY} + m_{\bfA} n\log(n)))$ flops.
The term $(k \wedge m) := \min(k,m)$ comes from the fact that if $k \geq m$ then we take $\bfOmega = \bfI_m$.

\paragraph{Sketched $HS(U_m,U_m)$ norm.}
We take $\bfGamma$ and $\bfSigma$ as fast to apply matrices, and $\bfOmega$ as a Gaussian, Rademacher or identity matrix.
We then evaluate the terms right to left in 
\[
    \bfLambda(\bfY_i \bfA_j)
    = \bfGamma \mathrm{vec}\big(
        (\bfSigma \bfU_m^* \bfA_j^*
        \bfY_i^*\bfU_m\bfOmega^*)^*
        \big).
\]
As a result, the total offline cost for $\|\bfE\|_{HS(U,U)}^{\bfPsi}$ is $\calO(p (k\wedge m) (T_{\bfY} + (k \wedge m) m_{\bfA}n\log(n)) )$ flops.
The term $(k \wedge m)$ comes from the fact that if $k \geq m$ then we take $\bfOmega = \bfSigma = \bfI_m$.
Not also that if $k \geq m^2$ then we can take all the matrices $\bfGamma$, $\bfOmega$ and $\bfSigma$ as identity matrices, which leads to $\bfLambda(\cdot) = \mathrm{vec}(\cdot)$ and to $\|\cdot\|_{HS(U_m, U_m)}^{\bfLambda} = \|\cdot\|_{HS(U_m, U_m)}$.

\subsection{Preconditioned Galerkin and residual}
\label{subsec:precond galerkin and residual}

\paragraph{Preconditioned Galerkin projection.}

Recall that the preconditioned Galerkin projection $\bfu_r(\bfxi) \in U_r$ defined in \eqref{equ:def precond gal} is written as $\bfu_r(\bfxi) = \bfU_r \bfa_r(\bfxi)$ where $\bfa_r \in \bbK^r$ is the solution to $\bfB_r(\bfxi) \bfa_r(\bfxi) = \bff_r(\bfxi)$ with $\bfB_r(\xi) = \bfU_r^* \bfR_U \bfP(\xi) \bfA(\xi) \bfU_r$ and $\bff_r(\xi) = \bfU_r^* \bfR_U \bfP(\xi) \bfb(\xi)$.
We can then write
\[
    \bfB_r(\bfxi) 
    = \sum_{i=1}^{p} \sum_{j=1}^{m_{\bfA}}
    \lambda_i(\bfxi) \theta_j^{\bfA}(\xi)
    \bfU_r^* \bfR_U \bfY_i \bfA_j \bfU_r,
    \quad
    \bff_r(\bfxi) 
    = \sum_{i=1}^{p} \sum_{j=1}^{m_{\bfb}}
    \lambda_i(\bfxi) \theta_j^{\bfb}(\xi)
    \bfU_r^* \bfR_U \bfY_i \bfb_j,
\]
where the coefficient vector $\bflambda(\xi) \in \bbK^p$ such that $\bfP(\xi) = \sum_{i=1}^p \lambda_i(\xi) \bfY_i$ is computed online by solving \eqref{equ:proj error A onto precond space} as described in \Cref{subsec:sketched measures of quality practical}.
Assuming that $pr \gtrsim \log(n)$, the total offline cost for computing the offline terms is then $\calO(p r (T_{\bfY} + n r m_{\bfA} + n m_{\bfb}))$.
The total online cost for summing the affine terms and solving the preconditioned Galerkin system is then $\calO( p r (r^2 + r m_{\bfA} + m_{\bfb}))$.

\paragraph{Projected preconditioned residual.}

Recall that the a posteriori error estimator from \Cref{prop:bound error by precond residual} is written as $\|\bfU_m^* \bfR_U \bfP(\bfxi)(\bfA(\bfxi) \bfu_r(\bfxi) - \bfb(\bfxi))\|_2$.
Writing $\bfu_r(\bfxi) = \bfU_r \bfa_r(\bfxi)$, the estimator can be written as $\|\bfB_{m,r}(\bfxi) \bfa_r(\xi) - \bff_m(\xi)\|_2$ where 
\[
    \bfB_{m,r}(\bfxi) 
    = \sum_{i=1}^{p} \sum_{j=1}^{m_{\bfA}}
    \lambda_i(\bfxi) \theta_j^{\bfA}(\xi)
    \bfU_m^* \bfR_U \bfY_i \bfA_j \bfU_r,
    \quad
    \bff_m(\bfxi) 
    = \sum_{i=1}^{p} \sum_{j=1}^{m_{\bfb}}
    \lambda_i(\bfxi) \theta_j^{\bfb}(\xi)
    \bfU_m^* \bfR_U \bfY_i \bfb_j,
\]
where the coefficient vector $\bflambda(\xi) \in \bbK^p$ such that $\bfP(\xi) = \sum_{i=1}^p \lambda_i(\xi) \bfY_i$ is computed online by solving \eqref{equ:proj error A onto precond space} as described in \Cref{subsec:sketched measures of quality practical}.
Assuming that $pm \gtrsim \log(n)$, the total offline cost for computing the offline terms is then $\calO(p m (T_{\bfY} + n r m_{\bfA} + n m_{\bfb}))$.
The total online cost for summing the affine terms is then $\calO( p m (r m_{\bfA} + m_{\bfb}))$.
Note that if $m$ is too large, then we can use a random sketching approach similar to \cite{balabanovRandomizedLinearAlgebra2019} to accurately and efficiently approximate the estimator on a finite parameter set, with a sketch size independent, or weakly dependent, on $m$.

\subsection{Greedy construction of the space of preconditioners}
\label{subsec:greedy construction preconditioner}

In this section, inspiring from \cite{zahmInterpolationInverseOperators2016}, we introduce a greedy algorithm to adaptively select the $\bfY_i = \bfA(\xi_i)^{-1}$, $i\geq 1$.
The error indicator the algorithm relies on has the form
\begin{equation}
\label{equ:error indicator for greedy}
    \Delta_i^{\bfTheta}(\xi) := 
    \min_{\bfP \in \spanv{\bfY_0, \cdots, \bfY_i}}
    \|\bfI - \bfP \bfA(\xi)\|_{HS}^{\bfTheta},
\end{equation}
with a sketched seminorm $\|\cdot\|^{\bfTheta}_{HS} \in \{\|\cdot\|_{HS(U, U)}^{\bfPsi}, \|\cdot\|_{HS(U, U_m)}^{\bfXi}, \|\cdot\|_{HS(U_m, U_m)}^{\bfLambda}\}$, or a weighted version $\|\cdot\|^{\bfTheta}_{HS} = ((\|\cdot\|_{HS(U_m, U_m)}^{\bfLambda})^2 + \tau^2 (\|\cdot\|_{HS(U, U_m)}^{\bfXi})^2)^{1/2}$ for some $\tau>0$ as mentioned in \Cref{rem:minimizing weighted sum}.
The algorithm we consider is then a classical greedy algorithm, as described in \Cref{algo:greedy preconditioner}.
Note that we take $\bfY_0 = \bfR_U^{-1}$, so that the space in which $\bfP$ is taken has dimension $p+1$ and so that iteration $0$ corresponds to the non-preconditioned setting.

\begin{algorithm}[H]
\caption{Greedy algorithm to build the space of preconditioners}
\label{algo:greedy preconditioner}
\begin{algorithmic}[1]
  \Require Orthonormal matrix $\bfU_m\in\bbK^{n \times m}$, error indicator $\Delta^{\bfTheta}$ from \eqref{equ:error indicator for greedy}, tolerance $\epsilon>0$, maximum dimension $p_{max}$.
  \Ensure Values $\bfY_1, \cdots, \bfY_p$ of $\bfA^{-1}(\xi)$ at interpolation points $\xi\in\{\xi_1, \cdots, \xi_p\}$.
  \State Set $i=0$ and $\bfY_0 = \bfR_U^{-1}$. 
  \While{$i\leq p_{max}$ and $\max_{\xi \in \calP} \Delta_i^{\bfTheta}(\xi) \geq \epsilon$}
    \State Set $i = i+1$.
    \State Set $\bfY_i = \bfA(\xi_i)^{-1}(\xi)$ with $\xi_i := \argmax_{\xi \in \calP}  \Delta_i^{\bfTheta}(\xi)$.
    \State Compute affine terms $\bfTheta(\bfY_i \bfA_j)$ for all $1\leq j\leq m_{\bfA}$ using \Cref{subsec:sketched measures of quality practical}.
    \State Compute affine terms for other sketched seminorms using \Cref{subsec:sketched measures of quality practical} and for preconditioned Galerkin and residual norm estimation using \Cref{subsec:precond galerkin and residual}.
  \EndWhile
  \end{algorithmic}
\end{algorithm}


\section{Numerical experiments}
\label{sec:numerics preconditioners}

\subsection{Setting}

In this section we consider the \emph{acoustic invisibility cloak} from 
\cite[Section 6.1]{balabanovRandomizedLinearAlgebra2021}, which is an acoustic wave scattering in 2D with a perfect scatterer covered in an invisibility cloak composed of layers of homogeneous isotropic materials.
Note that the implementation of the approach investigated in this work is available at \url{https://github.com/alexandre-pasco/rla4mor}.

\paragraph{Reduced spaces.}
The solutions of the PDE are approximated in a reduced space $U_r$ with dimension $r=50$.
The a posteriori error estimator is obtained by projecting the preconditioned residual onto $U_m$ with $m=250$.
The spaces $U_r$ and $U_m$ are obtained by performing the sketched proper orthogonal decomposition (POD) described in \cite[Section 5.2]{balabanovRandomizedLinearAlgebra2019}, truncating respectively the first $r$ and $m$ dominant modes, which implies that $U_r \subset U_m$.
The sketched POD is performed on snapshots computed for $2000$ random parameter values, using a P-SRHT sketch with $8192$ rows.

\paragraph{Sketch used for the measures of quality.}
The choice of the sketches for the HS norms follows the approach described in \Cref{subsec:sketched measures of quality practical}, taking $k=1024$ rows for all the random matrices considered.
For $\|\cdot\|_{HS(U, U)}^{\bfPsi}$, we take $\bfOmega$ as a Gaussian random matrix with $k$ rows, and we take $\bfGamma$ and $\bfSigma$ as P-SRHT matrices with $k$ rows.
For $\|\cdot\|_{HS(U, U_m)}^{\bfXi}$, we take $\bfOmega$ as the identity matrix $\bfI_m$, and we take $\bfGamma$ and $\bfSigma$ as P-SRHT matrices with $k$ rows.
For $\|\cdot\|_{HS(U_m, U_m)}^{\bfLambda}$, we take $\bfOmega$ and $\bfSigma$ identity matrices $\bfI_m$, and we take $\bfGamma$ as P-SRHT matrix with $k$ rows.

\paragraph{Greedy algorithms.}
We consider several instances of \Cref{algo:greedy preconditioner}.
Each instance build a space $\spanv{\bfY_0, \cdots, \bfY_{p_{max}}}$ by adaptively selecting the $\bfY_i = \bfA(\xi_i)^{-1}$ for $1\leq i\leq p_{max}=50$ according to the greedy criterion $\Delta^{\bfTheta}_i$ defined in \eqref{equ:error indicator for greedy}.
We compare greedy criterions using four different norms: $\|\cdot\|_{HS(U,U)}^{\bfPsi}$, $\|\cdot\|_{HS(U,U_m)}^{\bfXi}$, $\|\cdot\|_{HS(U_m,U_m)}^{\bfLambda}$ and $((\|\cdot\|_{HS(U_m,U_m)}^{\bfLambda})^2 + \tau^2 (\|\cdot\|_{HS(U,U_m)}^{\bfXi})^2)^{1/2}$ with $\tau^2 = 1/2$.
Due to computational time restrictions, the general purpose indicator $\|\cdot\|_{HS(U,U)}^{\bfPsi}$ was only computed until $p_{max} = 20$.
All the instances are performed on the same training set $\calP_{train}$ containing $\#\calP_{train} = 99000$ random parameter values.
Note that at iteration $0$ we have $\bfP = \bfY_0 = \bfR_U^{-1}$, which corresponds to the non-preconditioned setting.

\paragraph{Test quantities monitored.}
For all instances of the greedy algorithm, at each iteration, we monitor additional indicators evaluated on the same test set $\calP_{test}$ containing $\#\calP_{test} = 2000$ random parameter values, independent of $\calP_{train}$.
More precisely at step $j$, for every $\xi \in \calP_{test}$ we compute $\bfP(\xi)$ solution to \eqref{equ:error indicator for greedy}, the coefficients $\bfa_r(\xi) \in \bbK^r$ of the preconditioned Galerkin projection and we evaluate seven test quantities.
The quantiles of those test quantities are then summarized in the figures below.
The first three monitored quantities are related to the performance of the preconditioned Galerkin projection $\bfu_r$ defined in \Cref{prop:def precond gal} and the a posteriori error estimators.
We first monitor in \Cref{fig:errors vs greedy} an estimation of the normalized error,
\[
    \frac{
        \|\bfu(\xi) - \bfu_r(\xi)\|_U^{\bfPhi}
    }{
        \max_{\xi' \in \calP_{test}}
        \|\bfu(\xi')\|_U^{\bfPhi}
    },
\] 
where $\|\bfu(\xi) - \bfu_r(\xi)\|_U^{\bfPhi}$ and $\|\bfu(\xi')\|_U^{\bfPhi}$ are efficient sketched estimator of $\|\bfu(\xi) - \bfu_r(\xi)\|_U$ and $\|\bfu(\xi')\|_U^2$ respectively, using $\bfPhi = \bfPhi' \bfQ_U$ where $\bfPhi'$ is a P-SRHT matrix with $8192$ rows.
These sketched estimates are expected to be accurate with an error margin of $\pm 10 \%$.
We then monitor in \Cref{fig:quasi opti vs greedy} an estimation of the quasi-optimality of the preconditioned Galerkin projection,
\[
    \frac{\|\bfu(\xi) - \bfu_r(\xi)\|_U^{\bfPhi}}{\|\bfu(\xi) - \bfPi_{U_r} \bfu(\xi)\|_U} - 1.
\]
We then monitor in \Cref{fig:prnorms vs greedy} an estimation of the accuracy of the preconditioned residual-based error estimator from \Cref{prop:bound error by precond residual},
\[
    \max\big(
        \frac{
            \|\bfU_m^{\sfH} \bfR_U \bfP(\xi)(\bfA(\xi) \bfu_r(\xi) - \bfb(\xi))\|_2
        }{
            \|\bfu(\xi) - \bfu_r(\xi)\|_U^{\bfPhi}
        },
        \frac{
            \|\bfu(\xi) - \bfu_r(\xi)\|_U^{\bfPhi}
        }{
            \|\bfU_m^{\sfH} \bfR_U \bfP(\xi)(\bfA(\xi) \bfu_r(\xi) - \bfb(\xi))\|_2
        }
    \big) - 1.
\]
The next four monitored quantities are related to the measures of quality from \Cref{sec:measures of quality of a preconditioner}.
We monitor in \Cref{fig:uu vs greedy,fig:uum vs greedy,fig:umum vs greedy,fig:urur vs greedy} the measures of quality of the built preconditioner, namely
\[
    \|\bfE(\xi)\|_{HS(U,U)}^{\bfPsi},
    \quad \|\bfE(\xi)\|_{HS(U,U_m)}^{\bfXi},
    \quad \|\bfE(\xi)\|_{HS(U_m,U_m)}^{\bfLambda},
    \quad \|\bfE(\xi)\|_{HS(U_r, U_r)}.
\]
Recall that due to computational time restrictions,  $\|\bfE(\xi)\|_{HS(U,U)}^{\bfPsi}$ was only computed until $p_{max} = 20$.
Finally, the last quantity monitored in \Cref{fig:unstable vs greedy} is the frequency of the event $\|\bfE\|_{U_r, U_r} \geq 1$ on the test set, 
\[
    \frac{1}{\#\calP_{test}}
    \sum_{j=1}^{\#\calP_{test}}
    \mathbbm{1}_{\|\bfE\|_{U_r, U_r} \geq 1},
\]
under which the stability of the preconditioned Galerkin projection may not be satisfied.

\subsection{Detailed observations}
\label{subsec:detailed observations precond}

In this section we comment and analyze the numerical results obtained with the setting described in the previous section. 
We will in particular use the iteration $0$ as a baseline comparison as it corresponds to the non-preconditioned setting.

\begin{figure}[H]
    \centering
    \includegraphics[page=10, width=0.8\textwidth]{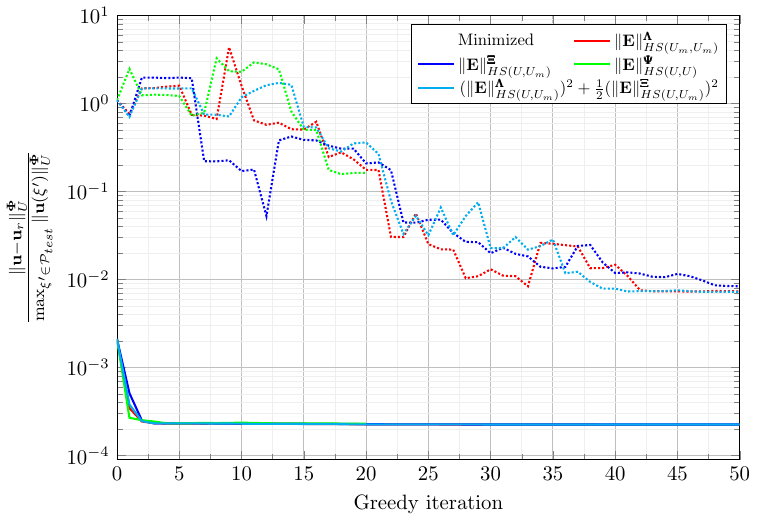}
    \caption[caption]{\footnotesize 
    Evolution of quantiles on the test set of the error along the greedy algorithms, using four different greedy criterions.
    The quantiles $90\%$ and $100\%$ are represented respectively by the continuous and dotted lines.
    }
    \label{fig:errors vs greedy}
\end{figure}

\begin{figure}[H]
    \centering
    \includegraphics[page=5, width=0.8\textwidth]{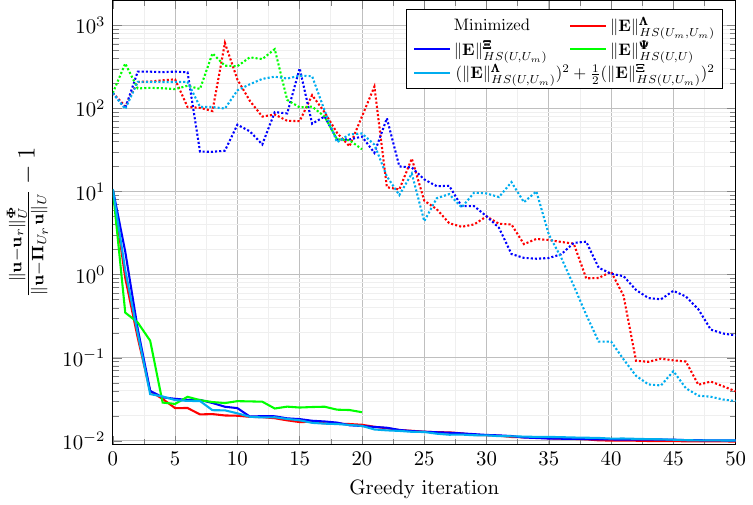}
    \caption[caption]{\footnotesize 
    Evolution of quantiles on the test set of the quasi-optimality gap of preconditioned Galerkin projection along the greedy algorithms, using four different greedy criterions.
    The quantiles $90\%$ and $100\%$ are represented respectively by the continuous and dotted lines.
    }
    \label{fig:quasi opti vs greedy}
\end{figure}

\begin{figure}[H]
    \centering
    \includegraphics[page=7, width=0.8\textwidth]{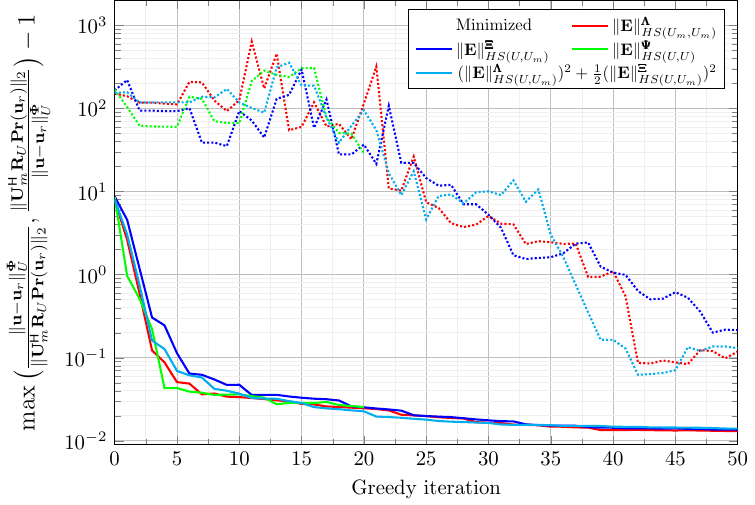}
    \caption[caption]{\footnotesize 
    Evolution of quantiles on the test set of the accuracy of the preconditioned residual-based error estimator along the greedy algorithms, using four different greedy criterions.
    The quantiles $90\%$ and $100\%$ are represented respectively by the continuous and dotted lines.
    }
    \label{fig:prnorms vs greedy}
\end{figure}

Firstly, we observe in \Cref{fig:errors vs greedy,fig:quasi opti vs greedy,fig:prnorms vs greedy} that for $90\%$ of the test samples, the four greedy criterions performed overall similarly, and it required few iterations to obtain significant improvements.
More precisely, with only $p=5$ iterations and for all greedy criterions, we make the following observations on $90\%$ of the test samples.
From \Cref{fig:errors vs greedy} the absolute error of the preconditioned Galerkin projection is smaller than $0.07$, compared to $0.62$ at $p=0$.
From \Cref{fig:quasi opti vs greedy} the quasi-optimality error of the preconditioned Galerkin projection is smaller than $1.5$, compared to $10.7$ at $p=0$.
From \Cref{fig:prnorms vs greedy} the accuracy of the preconditioned residual based error estimator is smaller that $0.12$, compared to $8.74$ at $p=0$.

Secondly, we observe in \Cref{fig:errors vs greedy,fig:quasi opti vs greedy,fig:prnorms vs greedy} that significant improvements were also obtained in the worst case scenarios, where the three model reduction oriented greedy criterions performed overall similarly, although it required more iterations.
More precisely, after $p=50$ iterations, and for the three model reduction oriented greedy criterions, we make the following observations on the worst cases of the test samples.
From \Cref{fig:errors vs greedy} the absolute error of the preconditioned Galerkin projection is smaller than $2.5$, compared to $319.9$ at $p=0$.
From \Cref{fig:quasi opti vs greedy} the quasi-optimality error of the preconditioned Galerkin projection is smaller than $0.2$, compared to $153.5$ at $p=0$.
This means that the preconditioned Galerkin projection is very close to the orthogonal projection, which can be observed in \Cref{fig:errors vs greedy} with the stagnation of the error during the last iterations.
From \Cref{fig:prnorms vs greedy} the accuracy of the preconditioned residual based error estimator is smaller that $0.22$, compared to $165.10$ at $p=0$.
This means that the preconditioned error estimator is very accurate.

\begin{figure}[H]
    \centering
    \includegraphics[page=8, width=0.8\textwidth]{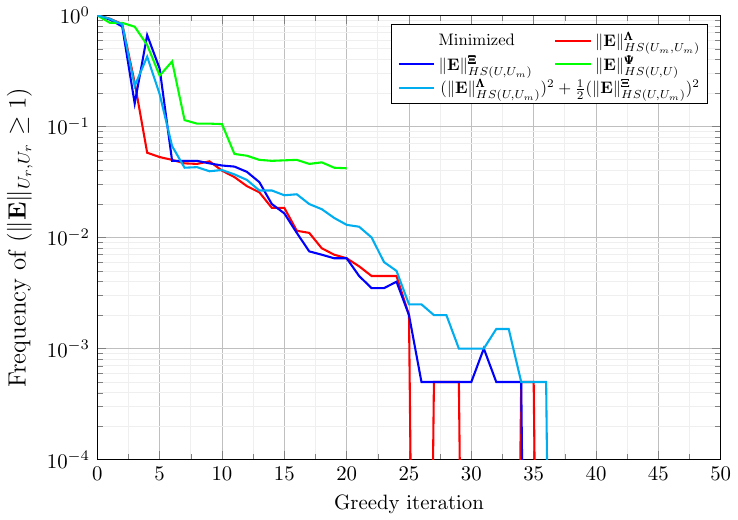}
    \caption[caption]{\footnotesize 
    Evolution of the frequency of the event $\|\bfI - \bfP\bfA\|_{U_r, U_r} \geq 1$ along the greedy algorithm, using four different greedy criterions.
    Under this event, the stability of the preconditioned Galerkin projection may not be satisfied.
    Frequency below $5\times 10^{-4} = 1/\#\calP_{test}$ means that the event never occurred.
    }
    \label{fig:unstable vs greedy}
\end{figure}

Thirdly, we observe in \Cref{fig:errors vs greedy} that during the first $15$ iterations, our preconditioning approach sometimes deteriorated the worst case error, for any greedy criterion.
This is most probably due to possible instability of the preconditioned Galerkin projection, which can occur when $\|\bfE\|_{U_r, U_r} \geq 1$.
We indeed observe in \Cref{fig:unstable vs greedy} that this situation remained plausible until $p=35$ iterations.
This problem could be circumvented by using the non-preconditioned Galerkin projection whenever $\|\bfE\|_{U_r, U_r} \geq 1$, or by updating $\bfP$ using an additional iterative minimization procedure on $\|\bfE\|_{U_r, U_r}$.

\begin{figure}[H]
    \centering
    \includegraphics[page=1, width=0.8\textwidth]{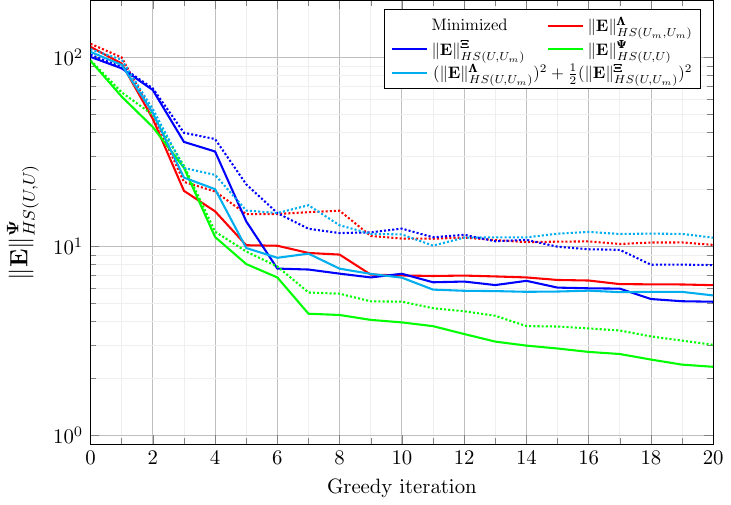}
    \caption[caption]{\footnotesize 
    Evolution of quantiles on the test set of $\|\bfI-\bfP(\xi)\bfA(\xi)\|_{HS(U,U)}^{\bfPsi}$ along the greedy algorithm, using four different greedy criterions.
    The quantiles $90\%$ and $100\%$ are represented respectively by the continuous and dotted lines.
    }
    \label{fig:uu vs greedy}
\end{figure}

\begin{figure}[H]
    \centering
    \includegraphics[page=2, width=0.8\textwidth]{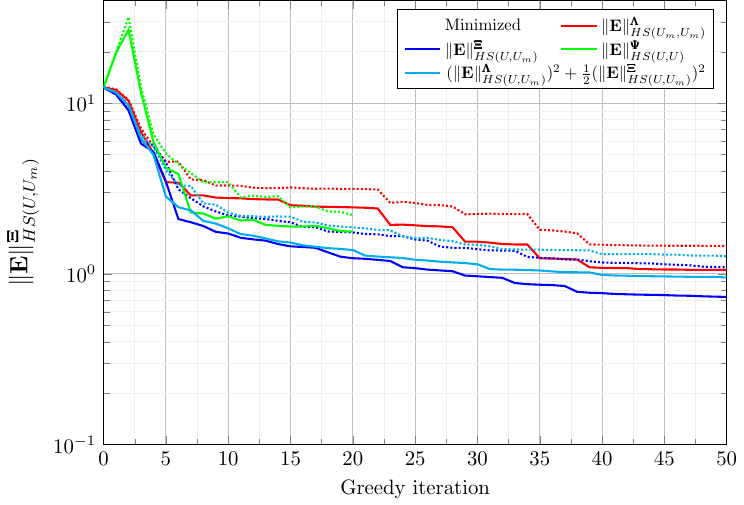}
    \caption[caption]{\footnotesize 
    Evolution of quantiles on the test set of $\|\bfI-\bfP(\xi)\bfA(\xi)\|_{HS(U,U_m)}^{\bfXi}$ along the greedy algorithm, using four different greedy criterions.
    The quantiles $90\%$ and $100\%$ are represented respectively by the continuous and dotted lines.
    }
    \label{fig:uum vs greedy}
\end{figure}

\begin{figure}[H]
    \centering
    \includegraphics[page=3, width=0.8\textwidth]{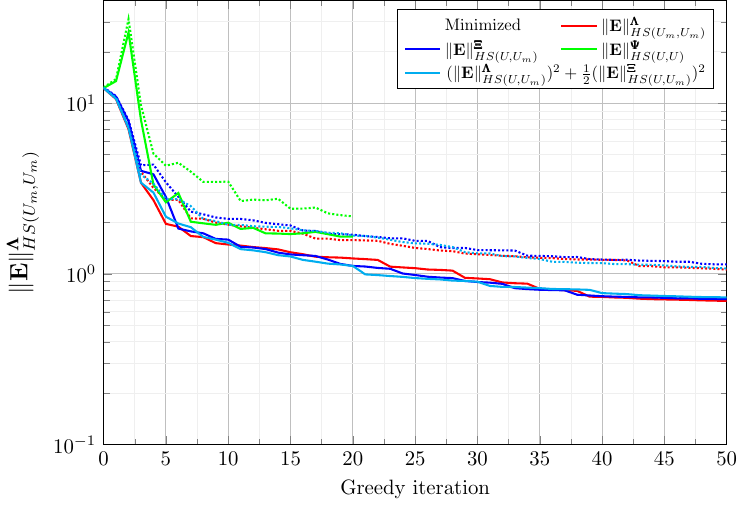}
    \caption[caption]{\footnotesize 
    Evolution of quantiles on the test set of $\|\bfI-\bfP(\xi)\bfA(\xi)\|_{HS(U_m,U_m)}^{\bfLambda}$ along the greedy algorithm, using four different greedy criterions.
    The quantiles $90\%$ and $100\%$ are represented respectively by the continuous and dotted lines.
    }
    \label{fig:umum vs greedy}
\end{figure}

\begin{figure}[H]
    \centering
    \includegraphics[page=4, width=0.8\textwidth]{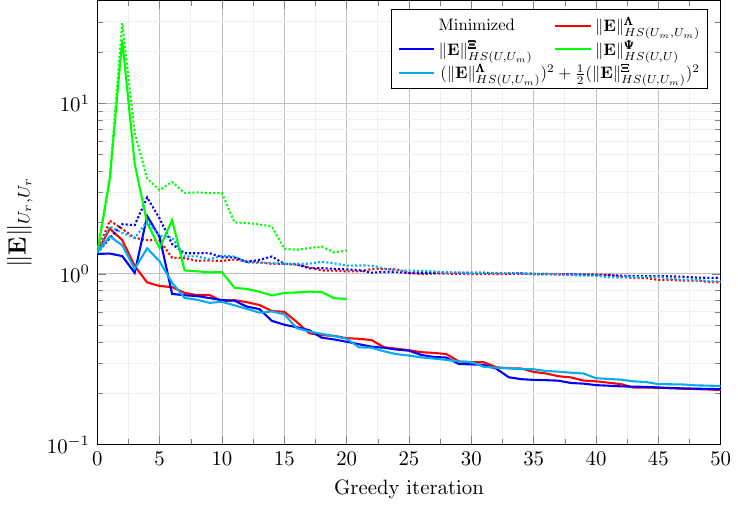}
    \caption[caption]{\footnotesize 
    Evolution of quantiles on the test set of $\|\bfI-\bfP(\xi)\bfA(\xi)\|_{U_r,U_r}$ along the greedy algorithm, using four different greedy criterions.
    The quantiles $90\%$ and $100\%$ are represented respectively by the continuous and dotted lines.
    }
    \label{fig:urur vs greedy}
\end{figure}

Fourthly, we observe in \Cref{fig:uu vs greedy,fig:uum vs greedy,fig:umum vs greedy,fig:urur vs greedy} that after a few iterations, there is an overall decrease of all the discrepancy measures.
The rate of decay may however depend on the greedy criterion.
For example in \Cref{fig:uu vs greedy} the model reduction oriented criterions yielded a slower decay of $\|\bfE\|_{HS(U,U)}^{\bfPsi}$ than the criterion $\|\bfE\|_{HS(U,U)}^{\bfPsi}$ after $p=6$ iterations.
This could be expected as $\|\cdot\|_{HS(U,U_m)}^{\bfXi}$ and $\|\cdot\|_{HS(U_m,U_m)}^{\bfLambda}$ are (much) weaker than $\|\cdot\|_{HS(U,U)}^{\bfPsi}$.
Also in \Cref{fig:uum vs greedy} the criterion $\|\bfE\|_{HS(U_m,U_m)}^{\bfLambda}$ performed the worst at minimizing $\|\bfE\|_{HS(U,U_m)}^{\bfXi}$ after $p=5$ iterations.
This could be expected as $\|\cdot\|_{HS(U_m,U_m)}^{\bfLambda}$ is weaker than $\|\cdot\|_{HS(U,U_m)}^{\bfXi}$.
On the other hand in \Cref{fig:umum vs greedy,fig:urur vs greedy} all the model reduction oriented criterions performed almost the same after $p=6$ iterations.
This is interesting, as we could expect the criterion $\|\bfE\|_{HS(U_m,U_m)}^{\bfLambda}$ to perform much better at minimizing $\|\bfE\|_{HS(U_m,U_m)}^{\bfLambda}$, but it turned out that the criterion $\|\bfE\|_{HS(U,U_m)}^{\bfXi}$ and the weighted sum performed just as well.

We thus conclude that, in our setting, the best choice of model reduction oriented greedy criterion seems to be $\|\bfE\|_{HS(U,U_m)}^{\bfXi}$ or $(\|\bfE\|_{HS(U_m,U_m)}^{\bfLambda})^2 + \tau^2(\|\bfE\|_{HS(U,U_m)}^{\bfXi})^2$ for some $\tau>0$, in our case $\tau = 1/\sqrt{2}$.

\section{Conclusion and perspectives}
\label{sec:conclusion preconditioners}

\subsection{Conclusion}

In the present work we proposed a new structured random sketching approach for embedding high-dimensional matrices given in implicit form, meaning that only input-output pairs of these matrices are accessible.
The new random embedding we proposed is constructed from three standard random matrices, which can be adapted to the matrix to embed, which makes the framework flexible.
We showed that it inherits the oblivious subspace embedding properties of the random matrices from which it is constructed, and provide  rigorous condition on their sizes for Gaussian and P-SRHT embeddings.
We also compared our approach with other existing methods, and showed that it is the best in terms of the number of input-output pairs required.
We then extended our approach to the case of Hilbert-Schmidt (HS) operators between finite dimensional Hilbert spaces.
Our approach allows us to approximate operators between high-dimensional spaces in a low-dimensional span of operators by solving a small linear least-squares problem.
Our analysis shows that the quasi-isometry of the embedding is satisfied with sketch sizes independent of the dimensions, which opens the way to an extension of our approach to the infinite dimensional setting.

As an application, we constructed preconditioners for high-dimensional linear equation by approximating the operator's inverse. 
We first showed that we can use various norms or seminorms on linear operators as measures of quality of a preconditioner, providing rigorous bounds on the quasi-optimality of the preconditioned Galerkin projection and accuracy of the residual-based a posteriori error estimator.
In particular, we introduced new seminorms tailored to the specific settings of projection-based model reduction, yielding sharper theoretical results.
Although all the theoretical results can be stated in terms of operator norms, those quantities are very challenging to minimize.

We mitigated this problem by considering HS norms, which bound the corresponding operator norms.
The main advantage is that such quantities define quadratic functions of the preconditioner, and since we take the latter within a low-dimensional space of operators, minimizing HS norms is equivalent to solving a linear least-squares problem.
However, there are two problems with the HS framework.
The first problem is that the operators are given in implicit form, making neither the least-squares system nor the residual tractably computable offline.
We circumvented this problem by using our new random sketching approach, leveraging the low-dimensionality of the space from which the preconditioners are taken.
The second problem is that the HS norm and the operator norm
can differ by a rather large factor, in general up to the square root of the full dimension.
We circumvented this problem by considering HS seminorms based on reduced spaces, which differ from the corresponding operator seminorms by a factor up to the square root of the reduced dimensions.

We also proposed a detailed practical approach for parameter dependent equations to obtain an efficient offline-online decomposition of the different parameter dependent quantities we introduced.
We illustrated our approach on an acoustic wave scattering as a numerical example.
We used a greedy algorithm to adaptively select interpolation points for the inverse of the full operator, which we used to define the space from which the preconditioner is taken.
We monitored, onto some random test set, quantiles of indicators such as the error, the quasi-optimality of the preconditioned Galerkin projection or the accuracy of the residual-based error estimator.
We then compared different choices of sketched HS norms as minimization criterions of the algorithms.
We observed that it required few greedy iterations to obtain good results on most of the test set.
Moreover, with enough greedy iterations we were able to obtain such results on the whole test set.
We also observed that, despite discrepancies during the very first iterations, all the greedy criterions essentially performed similarly.

\subsection{Perspectives}

Let us first mention two main perspectives to the current work concerning  operator sketching.
The first perspective is to extend our approach to low-rank tensor sketching, for which using structured embedding may be crucial.
The second perspective is to extend our approach to the infinite dimensional case, following for example \cite{boulleGeneralizationRandomizedSingular2022,boulleLearningEllipticPartial2023,perssonRandomizedNystromApproximation2025}.

Let us now mention three main perspectives to the current work concerning  preconditioning in a model order reduction setting.
The first perspective is to construct at the same time both the preconditioner space and the reduced solution space, as in \cite{zahmInterpolationInverseOperators2016}.
This raises practical questions, as whenever adding a direction to the reduced solution space, we would need to apply to it all the previously selected operators that span the current preconditioner space.
Following \Cref{rem:sketched mor} may allow an efficient approach.
The second perspective is to consider nonlinear approximation of the preconditioner, such as piecewise-affine or dictionary-based approximation.
The third perspective is to make more robust the overall preconditioned Galerkin projection estimator.
For example, whenever the HS seminorms are not sufficient to ensure good results, or even to ensure that the preconditioned reduced system is well posed, we could perform a further optimization procedure of the model reduction tailored operator seminorms. 
The sketched HS norms could then be seen as providing a good initialization to the problem of minimizing the corresponding operator norms.

Let us finally mention as a last perspective the efficient construction of preconditioners  for domain decomposition methods as in \cite{gosseletSimultaneousFETIBlock2015,spillaneAdaptiveMultiPreconditionedConjugate2016,bovetTwolevelAdaptationAdaptive2021}. Preconditioners for conjugate gradient methods are constructed as weighted sums of local preconditioners associated to each subdomain, and the choice of the weights can be very important to obtain fast convergence.  Then, our approach may be beneficial by allowing efficient optimization of these weights.

\appendix

\section{Preliminary results}
\label{apdx:preliminary}

In this section, we state three preliminary results for some Hilbert space $Z$ over a field $\bbK \in \{\bbR, \bbC\}$, equipped with an inner product $\innerp{\cdot}{\cdot}_Z$.
The results in this section are standard for the case $Z = \bbK^n$, see for example \cite{balabanovRandomizedLinearAlgebra2019} and supplementary materials, for which the proofs can be directly extended to a general space $Z$.
Although we provide the details for completeness, these results are not, strictly speaking, contributions of the present work.

The first preliminary result, \Cref{lem:approx norm implies approx inner} below, states that if a linear map approximates the norm on a subspace $V$ of $Z$, then it also approximates the inner product on $V$.
In the case $\bbK = \bbR$, it is a classical result that follows from parallelogram equalities.

\begin{lemma}
\label{lem:approx norm implies approx inner}
    Let $V$ be a subspace of $Z$.
    Let $\varepsilon>0$ and let $\bfTheta : Z \rightarrow \bbK^k$ be a linear map such that
    \[
        \forall \bfx \in V,\quad
        \left\vert \|\bfx\|_Z^2 - \|\bfTheta(\bfx)\|_2^2 \right\vert
        \leq \varepsilon \|\bfx\|_Z^2.
    \]
    Then, 
    \[
        \forall \bfx,\bfy \in V,\quad
        \left\vert \innerp{\bfx}{\bfy}_Z - \innerp{\bfTheta(\bfx)}{\bfTheta(\bfy)}_2 \right\vert
        \leq \varepsilon \|\bfx\|_Z \|\bfy\|_Z.
    \]
\end{lemma}
\begin{proof}
    Let $\bfx, \bfy \in V$.
    Firstly, it holds $\|\bfx + \bfy\|_Z^2 = \|\bfx\|_Z^2 + \|\bfy\|_Z^2 + 2 \mathrm{Re}(\innerp{\bfx}{\bfy}_Z)$, with $\mathrm{Re}(\cdot)$ the real part of a complex number.
    This yields the following parallelogram equalities
    \[
    \begin{aligned}
        \|\bfx + \bfy\|_Z^2 - \|\bfx - \bfy\|_Z^2
        & = 4 \mathrm{Re}(\innerp{\bfx}{\bfy}_Z),
        \\
        \|\bfx + \bfy\|_Z^2 + \|\bfx - \bfy\|_Z^2
        & = 2 \|\bfx\|_Z^2 + 2 \|\bfy\|_Z^2.
    \end{aligned}
    \]
    Let $\omega(\bfx,\bfy) := \innerp{\bfx}{\bfy}_Z - \innerp{\bfTheta(\bfx)}{ \bfTheta(\bfy)}_2 \in \bbK$.
    Using the first parallelogram identity with the linearity of $\bfTheta$ and the triangle inequality, we obtain 
    \[
        4 \left\vert \mathrm{Re}(\omega(\bfx,\bfy)) \right\vert
        \leq \left\vert
            \|\bfx + \bfy\|_Z^2 - \|\bfTheta(\bfx + \bfy)\|_2^2
            \right\vert 
            +
            \left\vert
            \|\bfx - \bfy\|_Z^2 - \|\bfTheta(\bfx - \bfy)\|_2^2
        \right\vert.
    \]
    Then, using the quasi-isometry property $\bfTheta$ with the facts that $\bfx + \bfy \in V$ and $\bfx - \bfy \in V$, then using the second parallelogram identity, we obtain 
    \[
        4 \left\vert \mathrm{Re}(\omega(\bfx,\bfy)) \right\vert
        \leq \varepsilon \|\bfx + \bfy\|_Z^2 + \|\bfx - \bfy\|_Z^2
        \leq 2 \varepsilon (\|\bfx\|_Z^2 + \|\bfy\|_Z^2).
    \]
    Using the above inequality for $\tilde \bfx = \|\bfx\|_Z^{-1}\bfx \in V$ and $\tilde \bfy = \|\bfy\|_Z^{-1}\bfy \in V$, and observing that $ \omega(\tilde \bfx, \tilde \bfy) \|\bfx\|_Z \|\bfy\|_Z  = \omega(\bfx,\bfy)$, we obtain
    \begin{equation}
    \label{equ:parallelogram bound real part}
        \left\vert \mathrm{Re}(\omega(\bfx,\bfy)) \right\vert
        \leq \varepsilon \|\bfx\|_Z \|\bfy\|_Z.
    \end{equation}
    If $\bbK = \bbR$, then this is the desired result since $\omega(\bfx,\bfy)\in\bbR$.
    Otherwise, if $\bbK = \bbC$, using the definition of $\omega$ with the Hermitian property of the inner product and the linearity of $\bfTheta$, we obtain
    \[
        |\omega(\bfx,\bfy)|^2
        = \omega(\bfx,\bfy)^* \omega(\bfx,\bfy)
        = \innerp{\omega(\bfx,\bfy)\bfx}{\bfy}_Z 
        - \innerp{\bfTheta(\omega(\bfx,\bfy) \bfx)}{ \bfTheta(\bfy)}_2
        = \omega(\omega(\bfx,\bfy)\bfx, \bfy),
    \]
    which implies that $\omega(\omega(\bfx,\bfy)\bfx, \bfy) \in \bbR$. Finally, using \eqref{equ:parallelogram bound real part} on $\omega(\bfx,\bfy) \bfx \in V$ and $\bfy \in V$, and using the property of the norm, we obtain
    \[
        |\omega(\bfx,\bfy)|^2
        = \mathrm{Re}(\omega(\omega(\bfx,\bfy)\bfx, \bfy))
        \leq \varepsilon \|\omega(\bfx,\bfy) \bfx\|_Z \|\bfy\|_Z
        = \varepsilon |\omega(\bfx,\bfy)| \|\bfx\|_Z \|\bfy\|_Z,
    \]
    which yields the desired result.
\end{proof}

The second preliminary result, \Cref{lem:decomposition on epsilon net} below, states that any element of a unit sphere can be decomposed as an infinite sum using any $\gamma$-net of this sphere, with geometric convergence of the sum.

\begin{lemma}
\label{lem:decomposition on epsilon net}
    Let $V$ be a subspace of $Z$.
    Let $\gamma \in (0,1)$ and let $\calN$ be a $\gamma$-net of $S := \{\bfx \in V : \|\bfx\|_Z=1\}$.
    Then for all $\bfx\in S$, there exists $(\bfy_i)_{i\geq 0} \subset \calN$ and $(\alpha_i)_{i\geq 0} \subset \bbR$ such that $0\leq \alpha_i \leq \gamma^i$ and 
    \begin{equation}
    \label{equ:decomposition on epsilon net}
        \forall m \in \bbN,
        \quad 
        \left\Vert \bfx - \sum_{i=0}^{m-1} \alpha_i \bfy_i \right\Vert_Z \leq \gamma^m.
    \end{equation}
\end{lemma}
\begin{proof}
    We prove this result by induction.
    Firstly, let $m=1$, by definition of $\calN$ we define $\bfy_0\in\calN$ such that $\|\bfx - \bfy_0\|_Z \leq \gamma$.
    Secondly, let $m=1$ and assume that we constructed $(\alpha_i, \bfy_i)_{0\leq m-1}$ such that $\| \bfx - \sum_{i=0}^{m-1} \alpha_i \bfy_i \|_Z \leq \gamma^m$.
    Let $\alpha_m := \|\bfx - \sum_{i=0}^{m-1} \alpha_i \bfy_i\|_Z$, which by assumption satisfies $\alpha_m \leq \gamma^m$.
    Then, consider $\bfz := \alpha_m^{-1} (\bfx - \sum_{i=0}^{m-1} \alpha_i \bfy_i)$.
    By definition, $\|\bfz\|_Z = 1$, and since $\bfx,\bfy_i \in V$ we have $\bfz \in V$, which yields $\bfz \in S$.
    Hence, by definition of $\calN$, there exists $\bfy_m \in \calN$ such that $\|\bfz - \bfy_m\|_Z \leq \gamma$.
    This yields
    \[
        \gamma^{m+1}
        \geq \alpha_m \gamma 
        \geq \alpha_m\|\bfz - \bfy_m\|_Z
        = \|\alpha_m \bfz - \alpha_m \bfy_m\|_Z
        = \left\Vert\bfx - \sum_{i=0}^m \alpha_i \bfy_i \right\Vert_Z.
    \]
    Hence, by induction, \eqref{equ:decomposition on epsilon net} holds.
\end{proof}

The third and last preliminary result, \Cref{lem:sketch on epsilon net} below, states that if a linear map approximates the inner product on a $\gamma$-net of a unit sphere, then it approximates the inner product on this sphere with an additional distortion factor smaller than $(1-\gamma)^{-2}$.

\begin{lemma}
\label{lem:sketch on epsilon net}
    Let $V$ be a $d$-dimensional subspace of $Z$.
    Let $\gamma \in (0,1)$ and let $\calN$ be a $\gamma$-net of $S := \{\bfx \in V : \|\bfx\|_Z=1\}$.
    Let $\eta>0$ and $\bfTheta : Z \rightarrow \bbR^k$ be a linear map such that 
    \[
        \forall \bfx,\bfy \in \calN,
        \quad 
        \left\vert
            \innerp{\bfx}{\bfy}_Z
            -\innerp{\bfTheta(\bfx)}{\bfTheta(\bfy)}_2
        \right\vert 
        \leq \eta.
    \]
    Then,
    \[
        \forall \bfx,\bfy \in V,
        \quad 
        \left\vert
            \innerp{\bfx}{\bfy}_Z
            -\innerp{\bfTheta(\bfx)}{\bfTheta(\bfy)}_2
        \right\vert 
        \leq \frac{\eta}{(1-\gamma)^2} \|\bfx\|_Z\|\bfy\|_Z.
    \]
\end{lemma}
\begin{proof}
    Let $\bfx \in V$.
    By \Cref{lem:decomposition on epsilon net} there exists $(\bfy_j)_{j\geq 0}\subset \calN$ and $(\alpha_j)_{j\geq 0} \subset \bbR$ such that $\bfx = \|\bfx\|_Z \sum_{i\geq 0} \alpha_j \bfy_j$ and $0 \leq \alpha_j \leq \gamma^j$.
    Firstly, using the assumptions on $\bfTheta$, then using the triangle inequality with $\alpha_j \geq 0$, and finally using $\alpha_j \leq \gamma^j$, we obtain
    \[
    \begin{aligned}
        \left\vert
            \|\bfTheta(\bfx)\|_2^2 
            - \|\bfx\|_Z^2 
        \right\vert
        &= \|\bfx\|_Z^2 \left\vert
            \sum_{j,l \geq 0} \left(
                \innerp{\bfTheta(\bfy_j)}{\bfTheta(\bfy_l)}_2 -  \innerp{\bfy_{j}}{\bfy_l}_Z
            \right) \alpha_j \alpha_l
        \right\vert 
        \\ 
        & \leq \|\bfx\|_Z^2 \eta \sum_{j,l\geq 0} \alpha_j \alpha_l
        = \|\bfx\|_Z^2 \eta \sum_{j\geq 0} \alpha_l)^2
        \leq \|\bfx\|_Z^2 \eta (1-\gamma)^{-2}.
    \end{aligned}
    \]
    As a result, $\bfTheta$ satisfies the assumptions of \Cref{lem:approx norm implies approx inner} with $\varepsilon=\eta (1-\gamma)^{-2}$, which yields the desired result.
\end{proof}


\section{Gaussian embedding on one side}
\label{apdx:gaussian ose one side}

In this section, we prove \Cref{prop:ose with dim free gaussian}, which states a non-asymptotic result for the oblivious $\bbK^{q\times p} \rightarrow \bbK^{q \times k}$ subspace embedding property of $\bfX \mapsto \bfX \bfOmega^*$ when $\bfOmega \in \bbR^{k\times p}$ is drawn from the rescaled Gaussian distribution.
Note that one should be able to show similar results when the entries of $\bfOmega$ are independently drawn from a general sub-Gaussian distribution, such as the Rademacher distribution, by using the \emph{Hanson-Wright} inequality which generalizes \Cref{lem:norm Aw is sub exp} below, see for example \cite[Section 6.2]{vershyninHighDimensionalProbabilityIntroduction2018}, up to modification in the constants.

\begin{lemma}
\label{lem:norm Aw is sub exp}
    Let $\bfX \in \bbR^{q\times p}$ and $\bfomega \sim \calN(0,\bfI_p)$.
    Then $\|\bfX \bfomega\|_2^2$ is sub-exponential with parameters $(2\|\bfX\|_{\sigma,4}^2, 4 \sigma_1(\bfX)^2)$, where $\|\bfX\|_{\sigma,4}$ is the Schatten $4$-norm of $\bfX$.
\end{lemma}
\begin{proof}
    Without any loss of generality, we assume that $q \leq p$.
    Let $\bfX = \bfU \bfS \bfV^T$ be the SVD of $\bfX$, with $\bfV \in \bbR^{q \times p}$ and $\bfU \in \bbR^{q\times q}$ such that $\bfV^T \bfV = \bfI_q$ and $\bfU^T \bfU = \bfI_q$, and with $\bfS=\mathrm{diag}(\sigma_1, \cdots, \sigma_q)$ such that $\sigma_1 \geq \cdots \geq \sigma_q \geq 0$.
    We have
    \[
        \|\bfX\bfomega\|_2^2 
        = \|\bfS \bfz\|_2^2
        = \sum_{i=1}^q \sigma_i^2 z_i^2,
    \]
    where $\bfz := \bfV^T \bfomega \in \bbR^q$ is a Gaussian vector with zero mean and covariance matrix
    \[
        \Expe{\bfz \bfz^T}
        = \Expe{\bfV^T \bfomega \bfomega^T \bfV}
        = \bfV^T \Expe{ \bfomega \bfomega^T }\bfV
        = \bfV^T \bfI_p \bfV
        = \bfI_q.
    \]
    Hence, $(z_i)_{1\leq i\leq q}$ are i.i.d with $z_i \sim \calN(0,1)$.
    We can now use \cite[Section 2.1]{wainwrightHighDimensionalStatisticsNonAsymptotic2019} to conclude the proof.
    Indeed, $(z_i^2)_{1\leq i\leq q}$ are i.i.d. sub-exponential with parameters $(2, 4)$.
    Also, $(\sigma_i^2 z_i^2)_{1\leq i\leq q}$ are i.i.d. sub-exponential with respective parameters $(2\sigma_i^2, 4\sigma_i^2)$.
    Finally, $\|\bfX\bfomega\|_2^2 = \sum_{i=1}^q \sigma_i^2 z_i^2$ is sub-exponential with parameters $(2 \sqrt{\sum_{i=1}^q \sigma_i^4}, 4 \max_{1\leq i\leq q} \sigma_i)$, where $\sqrt{\sum_{i=1}^q \sigma_i^4} = \|\bfX\|_{\sigma,4}^2$ and $\sigma_1 \geq \sigma_i$.
\end{proof}

From the above \Cref{lem:norm Aw is sub exp}, we can then show a concentration inequality for $\|\bfX\bfOmega^*\|_F^2$, as stated in \Cref{lem:norm Aw concentration} below.

\begin{lemma}
\label{lem:norm Aw concentration}
    Let $\bfX \in \bbR^{q\times p}$, let $k\in\bbN^*$,
    and let $\bfOmega = k^{-1/2}(\bfomega^{(1)}, \cdots, \bfomega^{(k)})^* \in \bbR^{k\times p}$ with $\bfomega^{(1)}, \cdots, \bfomega^{(k)}$ i.i.d. random vectors with $\bfomega^{(1)} \sim \calN(0,\bfI_p)$.
    Then for all $\varepsilon\in(0,1)$,
    \[
    \begin{aligned}
        \Proba{\left\vert\|\bfX \bfOmega^*\|_F^2 - \|\bfX\|_F^2 \right\vert > \varepsilon \|\bfX\|_F^2}
        &\leq 2 \exp\left(
            -\frac{k}{8} 
            \min\left\{
                \frac{\varepsilon^2 \|\bfX\|_F^4}{\|\bfX\|_{\sigma,4}^4},
                \frac{\varepsilon \|\bfX\|_F^2}{\sigma_1(\bfX)^2}
            \right\}
        \right) 
        \\
        &\leq 2 \exp\left(
            -\frac{k \varepsilon^2}{8}
        \right).
    \end{aligned}
    \]  
\end{lemma}
\begin{proof}
    Let $Z := \|\bfX \bfOmega^*\|_F^2 = \frac{1}{k}\sum_{i=1}^{k}\|\bfX \bfomega^{(i)}\|^2_2$, which satisfies $\Expe{Z} = \|\bfX\|_F^2$ by definition of $\bfomega^{(i)}$.
    Using \Cref{lem:norm Aw is sub exp}, the random variables $(\|\bfX \bfomega^{(i)}\|^2_2)_{1\leq i\leq k}$ are i.i.d. sub-exponential with parameters $(2\|\bfX\|_{\sigma,4}^2, 4\sigma_1(\bfX)^2)$.
    Hence, $Z$ is sub-exponential with parameters $(2 k^{-1/2} \|\bfX\|_{\sigma,4}^2, 4 k^{-1} \sigma_1(\bfX)^2)$, which implies
    \[
        \forall \eta >0, \quad
        \Proba{|Z - \Expe{Z}| \geq \eta}
        \leq 2 \exp\left(
            -\frac{k}{8} 
            \min\left\{
                \frac{\eta^2}{\|\bfX\|_{\sigma,4}^4},
                \frac{\eta}{\sigma_1(\bfX)^2}
            \right\}
        \right).
    \]
    Then, for all $\varepsilon>0$, applying the above inequality to $\eta = \|\bfX\|_F^2 \varepsilon > 0$ yields the first desired inequality.
    Finally, the right-hand side in the above inequality is a decreasing function of both $\|\bfsigma(\bfX)\|_4^4$ and $\sigma_1(\bfX)^2$.
    Hence, using that $\sigma_1(\bfX)^2 \leq \|\bfX\|_F^2$ and 
    \[
        \|\bfX\|_{\sigma,4}^4 = \sum_{i=1}^q \sigma_i^4  \leq \left( \sum_{i=1}^q \sigma_i^2 \right)^2 = \|\bfX\|_F^4,
    \]
    we obtain $\varepsilon^2 \|\bfX\|_F^4 / \|\bfX\|_{\sigma,4}^4 \leq \varepsilon^2$ and $\varepsilon \|\bfX\|_F^2 /\sigma_1(\bfX)^2 \leq \varepsilon \leq \varepsilon^2$ for all $\varepsilon\in(0,1)$, which yields the second desired inequality.
\end{proof}

It is worth noting that in the first inequality in \Cref{lem:norm Aw concentration}, the concentration inequality depends on the spectrum of the matrix to sketch.
In particular, the fastest concentration is obtained for a matrix with flat spectrum.
For example take $\bfX = \bfI_p$, then the upper bound in \Cref{lem:norm Aw concentration} is $2 \exp(-kp \varepsilon^2/8)$ for all $\varepsilon>0$ small enough.
This can be expected since in that case, $\|\bfX \bfOmega^*\|_F^2 = \| \bfOmega^*\|_F^2 = \|\mathrm{vec}(\bfOmega)\|_2^2$, where $\mathrm{vec}(\bfOmega)$ is a random vector with $kp$ i.i.d. $\calN(0, k^{-1/2})$ entries.

On the other hand, the slowest concentration is obtained for rank-one matrices, for which the upper bound in \Cref{lem:norm Aw concentration} is $ 2\exp(-k \varepsilon^2/8)$, which is essentially the same as the classical concentration for sketching vectors with a Gaussian sketch.
This can be expected since in that case, writing $\bfX = \bfu \bfv^*$ with $\|\bfu\|_2=\|\bfv\|_2=1$ yields $\|\bfX \bfOmega^*\|_F^2 = \| \bfOmega \bfu\|_2^2$.

Finally, using the results from \Cref{apdx:preliminary}, we can proof the main result of this section, which is \Cref{prop:ose with dim free gaussian} below.

\begin{proposition}
\label{prop:ose with dim free gaussian}
    Let $\bfOmega \in \bbR^{k\times p}$ with i.i.d. $\calN(0, k^{-1/2})$ entries, with
    \[
        k \geq 10.5 \varepsilon^{-2}( 6.9 \eta_{\bbK} d + \log(2/\delta)),
        \quad \eta_{\bbR} = 1,
        \quad \eta_{\bbC} = 2.
    \]
    Then $\bfX \mapsto \bfX \bfOmega^*$ is an $(\varepsilon, \delta, d)$ oblivious $\bbK^{q\times p} \rightarrow \bbK^{q\times k}$ subspace embedding.
\end{proposition}
\begin{proof}
    Let us first assume that $\bbK = \bbR$.
    Let $\gamma\in(0,1)$, then from \cite[Lemma 2.4]{bourgainApproximationZonoidsZonotopes1989} there exists a $\gamma$-net $\calN \subset S:=\{\bfX \in V : \|\bfX\|_F=1\}$ such that $\#\calN \leq (1+2/\gamma)^d$.
    Let $\eta \in (0,1)$ and let $\bfOmega \in \bbR^{k\times p}$ with i.i.d Gaussian entries, with zero mean and variance $1/k$, with $k$ such that
    \[
        k 
        \geq 8 \eta^2 \left(
            2d \log(1 + 2/\gamma) + \log(2/\delta)
        \right)
        \geq 8 \eta^2 \left(
            \log(\#\calN^2) + \log(2/\delta)
        \right).
    \]
    Then, \Cref{lem:norm Aw concentration} yields that for any $\bfX \in V$, 
    \[  
        \Proba{\left\vert\|\bfX \bfOmega^*\|_F^2 - \|\bfX\|_F^2 \right\vert > \eta \|\bfX\|_F^2}
        \leq \delta (\#\calN)^{-2}.
    \]
    Now, consider the set 
    $
        \widetilde\calN := 
        \{\bfx + \bfy : \bfx,\bfy \in \calN\} 
        \cup \{\bfx - \bfy : \bfx,\bfy \in \calN\},
    $
    which contains less than $(\#\calN)(\#\calN-1) +1 \leq (\#\calN)^2$ distinct directions. 
    Since $\calN\subset V$, we have $\widetilde\calN \subset V$.
    Then, using the above inequality, and using a union bound argument, we obtain that 
    \begin{equation}
    \label{equ:isometrie on plus minus net}
        \forall \bfX \in \widetilde \calN,\quad
        \left\vert\|\bfX\bfOmega^*\|_F^2 - \|\bfX \|_F^2 \right\vert \leq \eta \|\bfX\|_F^2
    \end{equation}
    holds with probability larger than $1-\delta$.
    Moreover, if \eqref{equ:isometrie on plus minus net} holds, then using the parallelogram equalities as in \Cref{lem:approx norm implies approx inner}, we obtain 
    \[
        \forall \bfX,\bfY \in \calN,\quad
        \left\vert
            \innerp{\bfX}{\bfY}_F - \innerp{\bfX \bfOmega^*}{\bfY \bfOmega^*}_F 
        \right\vert 
        \leq \eta \|\bfX\|_F  \|\bfY\|_F = \eta,
    \]
    which holds with probability larger than $1-\delta$.
    Now, if \eqref{equ:isometrie on plus minus net} holds, then by applying \Cref{lem:sketch on epsilon net} with $\bfTheta: \bfX \mapsto \mathrm{vec}(\bfX \bfOmega^*)$, we obtain that
    \[
        \forall \bfX,\bfY \in V,
        \quad
        \left\vert
            \innerp{\bfX}{\bfY}_F - \innerp{\bfX \bfOmega^*}{\bfY \bfOmega^*}_F 
        \right\vert
        \leq \eta (1-\gamma)^{-2} \|\bfX\|_F  \|\bfY\|_F,
    \]
    holds with probability larger than $1-\delta$.
    As a result, for all $\gamma \in (0,1)$ and all $\varepsilon \in (0,1)$, applying the previous result with $\eta = (1-\gamma)^{2} \varepsilon \in (0,1)$ yields that if
    \[
        k \geq 8 \varepsilon^{-2} \frac{2d \log(1+2/\gamma) + \log(2/\delta)}{(1-\gamma)^4},
    \]
    then
    \[
        \Proba{
        \forall \bfX,\bfY \in V,
        \quad
        \left\vert
            \innerp{\bfX}{\bfY}_F - \innerp{\bfX \bfOmega^*}{\bfY \bfOmega^*}_F 
        \right\vert
        \leq \varepsilon \|\bfX\|_F  \|\bfY\|_F
        } \geq 1-\delta.
    \]
    Now, when $2d \log(1+2/\gamma) \ll \log(2/\delta)$, the value $\gamma$ that minimizes the lower bound on $k$ is close to $\argmin_{x \in (0,1)}\log(1+2/x)(1-x)^{-4}$, which is approximately $\hat\gamma:=0.0656$.
    Replacing $\gamma$ with $\hat \gamma$ in the lower bound for $k$ yields the desired result for $\bbK = \bbR$.

    Finally, for the case $\bbK = \bbC$, we can use the same approach as in the supplementary material of \cite{balabanovRandomizedLinearAlgebra2019}.
    By using the triangle inequality with the fact that $\bfOmega$ is a real matrix, we obtain
    \[
        \left\vert
            \|\bfX\|_F^2 - \|\bfX\bfOmega^*\|_F^2
        \right\vert
        \leq
        \left\vert
            \|\mathrm{Re}(\bfX)\|_F^2 - \|\mathrm{Re}(\bfX)\bfOmega^*\|_F^2
        \right\vert 
        + \left\vert
            \|\mathrm{Im}(\bfX)\|_F^2 - \|\mathrm{Im}(\bfX)\bfOmega^*\|_F^2
        \right\vert. 
    \]
    Then, let $(\bfV_j)_{1\leq j\leq d}$ be a basis of $V$ and let $W$ be the space defined by
    \[
        W := \spanv{\mathrm{Re}(\bfV_1), \cdots, \mathrm{Re}(\bfV_d), \mathrm{Im}(\bfV_1), \cdots, \mathrm{Im}(\bfV_d)},
    \]
    which is a $2d$-dimensional subspace of $\bbR^{q\times p}$, and satisfies $\mathrm{Re}(\bfX)\in W$ and $\mathrm{Im}(\bfX) \in W$,  for all $\bfX \in V$.
    Using the previous inequality and the fact that $\bfX \mapsto \bfX \bfOmega^*$ is an $(\varepsilon, \delta, \eta_{\bbC}d)$ oblivious subspace embedding, where $\eta_{\bbC}=2$, we obtain
    \[
        \Proba{
            \forall \bfX \in V,~
            \left\vert
                \|\bfX\|_F^2 - \|\bfX\bfOmega^*\|_F^2
            \right\vert
            \leq \varepsilon \|\bfX\|_F^2
        } \geq 1-\delta.
    \]
    Finally, applying \Cref{lem:approx norm implies approx inner} to $\bfTheta: \bfX \mapsto \mathrm{vec}(\bfX\bfOmega^*)$ yields the desired result.
\end{proof}


\section{General embedding on one side}
\label{apdx:general ose one side}

In this section, we prove \Cref{prop:ose union bound} and \Cref{prop:ose gaussian trick}, which both state that if $\bfOmega$ is an oblivious $\bbK^p \rightarrow \bbK^k$ subspace embedding, then $\bfX \mapsto \bfX \bfOmega^*$ is an oblivious $\bbK^{q\times p} \rightarrow \bbK^{q\times k}$ subspace embedding, with slightly deteriorated properties.
The difference between \Cref{prop:ose union bound} and \Cref{prop:ose gaussian trick} lies in the deterioration of the embedding properties.
Indeed, the former introduces a factor $q$ in the probability of failure, while  the latter introduces factors in the accuracy and the probability of failure.

The most direct approach to obtain an oblivious embedding property on $\bfX \mapsto \bfX\bfOmega^*$ is to use a union bound argument on the rows of $\bfX$, as used in \Cref{prop:ose union bound} below.

\begin{proposition}
\label{prop:ose union bound}
    Let $\bfOmega \in \bbK^{k\times p}$ be an $(\varepsilon, \delta, d)$ oblivious $\bbK^p \rightarrow \bbK^k$ subspace embedding.
    Then, $\bfOmega$ is an $(\varepsilon, q\delta, d)$ oblivious $\bbK^{q\times p} \rightarrow \bbK^{q \times k}$ subspace embedding. 
\end{proposition}
\begin{proof}
    Consider for all $1\leq i\leq q$ the space 
    \[
        V_i := \{\bfX^* \bfe_i: \bfX\in V\} \subset \bbK^p,
    \]
    where $\bfe_i\in\bbK^q$ denotes the $i$-th column of the $q\times q$ identity matrix.
    Since $V_i$ is a $d$-dimensional subspace of $\bbK^p$, using the definition of $\bfOmega$ and a union bound argument, we obtain
    \[
        \Proba{
            \forall i\in\{1, \cdots, q\},~
            \forall \bfx\in V_i,~
            \big|\|\bfx\|_2^2 - \|\bfOmega \bfx\|_2^2\big|
            \leq \varepsilon \|\bfx\|_2^2
        } \geq 1- q\delta,
    \]
    which, by definition of $V_i$, is equivalent to 
    \[
        \Proba{
            \forall i\in\{1, \cdots, q\},~
            \forall \bfX \in V,~
            \big|\|\bfX^* \bfe_i\|_2^2 - \|\bfOmega \bfX^* \bfe_i\|_2^2\big|
            \leq \varepsilon \|\bfX \bfe_i\|_2^2
        } \geq 1- q\delta.
    \]
    Then, using the above equation, using
    $
        \|\bfX \|_F^2 - \|\bfX \bfOmega^* \|_F^2 
        = \sum_{i=1}^q (\|\bfX^* \bfe_i\|_2^2 - \|\bfOmega \bfX^* \bfe_i\|_2^2)
    $
    and using the triangle inequality, we obtain,
    \[
        \Proba{
            \forall \bfX \in V,~
            \big|\|\bfX \|_F^2 - \|\bfX\bfOmega^* \|_F^2\big|
            \leq \varepsilon \|\bfX \|_F^2
        } \geq 1- q\delta.
    \]
    Finally, applying \Cref{lem:approx norm implies approx inner} to $\bfTheta: \bfX \mapsto \mathrm{vec}(\bfX\bfOmega^*)$ yields the desired result.
\end{proof}

The problem with the above \Cref{prop:ose union bound} is that it depends on the dimension of the matrices to be sketched, via the probability of failure.
Although this dependency is expected to be logarithmic, thus negligible in many practical applications, it can still be problematic, when for example trying to generalize to the infinite-dimensional setting.

This problem can be addressed with an intermediate oblivious $\bbK^{q\times p} \rightarrow \bbK^{k_0 \times p}$ subspace embedding $\bfX \mapsto \bfOmega_0 \bfX$ with $\bfOmega_0 \in \bbK^{k_0 \times q}$, and essentially write
\[
    \|\bfX\|_F^2 
    \approx \|\bfOmega_0 \bfX\|_F^2
    \approx \|\bfOmega_0 \bfX \bfOmega^*\|_F^2
    \approx \|\bfX \bfOmega^*\|_F^2,
\]
using the union bound argument from \Cref{prop:ose union bound} to write $\|\bfOmega_0 \bfX\|_F^2 \approx \|\bfOmega_0 \bfX \bfOmega^*\|_F^2$ where the size of $\bfOmega$ depends on $k_0$ but not on $q$.
This is detailed in \Cref{prop:ose gaussian trick forward} below.

\begin{lemma}
\label{prop:ose gaussian trick forward}
    Let $\bfOmega \in \bbK^{k\times p}$ be an $(\varepsilon', \delta', d)$ oblivious $\bbK^p \rightarrow \bbK^k$ subspace embedding.
    Let $\varepsilon_0,\delta_0\in [0,1)$, and assume that there exists an $(\varepsilon_0, \delta_0, d)$ oblivious $\bbK^{q\times p} \rightarrow \bbK^{k_0 \times p}$ subspace embedding for some $k_0 >0$, independent of $\bfOmega$.
    Then, $\bfOmega$ is an $(\varepsilon, \delta, d)$ oblivious $\bbK^{q\times p} \rightarrow \bbK^{q \times k}$ subspace embedding, with $\varepsilon := \frac{1}{1-\varepsilon_0}(2\varepsilon_0 + \varepsilon' (1+\varepsilon_0))$ and $\delta:= k_0\delta' + 2\delta_0$.
\end{lemma}
\begin{proof}
    Let $\bfOmega_0\in\bbK^{k_0 \times q}$ be such that $\bfX \mapsto \bfOmega_0 \bfX$ is an $(\varepsilon_0, \delta_0, d)$ oblivious $\bbK^{q\times p} \rightarrow \bbK^{k_0 \times p}$ subspace embedding.
    Define $V(\bfOmega_0) := \{\bfOmega_0\bfX: \bfX \in V\}$ which is a $d$ dimensional subspace of $\bbK^{k_0 \times p}$.
    Using the independence of $\bfOmega$ and $\bfOmega_0$, and using \Cref{prop:ose union bound}, we obtain
    \[
    \begin{aligned}
        & \Proba{
            \forall \bfY \in V(\bfOmega_0),~
            \big|\|\bfY \|_F^2 - \|\bfY\bfOmega^* \|_F^2\big|
            \leq \varepsilon \|\bfY \|_F^2
        }
        \\ 
        & \quad = \Expe[\bfOmega_0]{
            \Proba[\bfOmega | \bfOmega_0]{
            \forall \bfY \in V(\bfOmega_0),~
            \big|\|\bfY \|_F^2 - \|\bfY\bfOmega^* \|_F^2\big|
            \leq \varepsilon \|\bfY \|_F^2
            }
        }
        \geq 1 - k_0 \delta'.
    \end{aligned}
    \]
    Secondly, using the embedding property of $\bfOmega_0$, we obtain
    \[
        \Proba{
            \forall \bfX \in V,~
            \big|\|\bfX \|_F^2 - \|\bfOmega_0\bfX \|_F^2\big|
            \leq \varepsilon_0 \|\bfX \|_F^2
        } \geq 1-\delta_0.
    \]
    Now, define $W(\bfOmega) := \{\bfX \bfOmega^* : \bfX \in V\}$ and $\tilde W(\bfOmega) := \{ \bfX (\bfOmega^* , \mathbf{0}_{p-k}) : \bfX \in V\}$, which are $d$ dimensional subspaces of $\bbK^{q \times k}$ and $\bbK^{q \times p}$ respectively.
    Using the independence of $\bfOmega$ and $\bfOmega_0$, and using the embedding property of $\bfOmega_0$, we obtain
    \[
    \begin{aligned}
        & \Proba{ 
            \forall \bfY \in \tilde W(\bfOmega),~
            \big|\|\bfY \|_F^2 - \|\bfOmega_0\bfY \|_F^2\big|
            \leq \varepsilon_0 \|\bfY \|_F^2
        }
        \\ 
        & \quad = \Expe[\bfOmega]{
            \Proba[\bfOmega_0 | \bfOmega]{
            \forall \bfY \in \tilde W(\bfOmega),~
            \big|\|\bfY \|_F^2 - \|\bfOmega_0\bfY \|_F^2\big|
            \leq \varepsilon_0 \|\bfY \|_F^2
            }
        }
        \geq 1 - \delta_0.
    \end{aligned}
    \]
    Moreover since $\|\bfOmega_0 \bfX (\bfOmega^* , \mathbf{0}_{p-k})\|_F^2 = \|\bfOmega_0 \bfX \bfOmega^*\|_F^2$, the above inequality also holds if we replace $\tilde W(\bfOmega)$ with $W(\bfOmega)$.
    As a result, by bounding the probability of failure of one of the three subspace embedding events, we obtain that these three events hold simultaneously with probability at least $1 - \delta$ with $\delta:=k_0\delta' + 2\delta_0$.
    Under these events, it holds
    \[
    \begin{gathered}
        \frac{1}{1+\varepsilon_0}
        \|\bfOmega_0 \bfX \bfOmega^*\|_F^2
        \leq \|\bfX \bfOmega^*\|_F^2
        \leq \frac{1}{1-\varepsilon_0}
        \|\bfOmega_0 \bfX \bfOmega^*\|_F^2,
        \\
        \frac{1}{1+\varepsilon_0} (1-\varepsilon')
        \|\bfOmega_0 \bfX\|_F^2
        \leq \|\bfX \bfOmega^*\|_F^2
        \leq \frac{1}{1-\varepsilon_0} (1+\varepsilon')
        \|\bfOmega_0 \bfX\|_F^2,
        \\
        \frac{1-\varepsilon_0}{1+\varepsilon_0} (1-\varepsilon')
        \| \bfX \|_F^2
        \leq \|\bfX\bfOmega^*\|_F^2
        \leq \frac{1+\varepsilon_0}{1-\varepsilon_0} (1+\varepsilon')
        \| \bfX\|_F^2,
        \\
        (1-\varepsilon)
        \| \bfX \|_F^2
        \leq \|\bfX\bfOmega^*\|_F^2
        \leq  (1+\varepsilon)
        \| \bfX\|_F^2,
     \end{gathered}  
    \]
    with $\varepsilon := \frac{1}{1-\varepsilon_0}(2\varepsilon_0 + \varepsilon' (1+\varepsilon_0))$.
    Finally, applying \Cref{lem:approx norm implies approx inner} to $\bfTheta: \bfX \mapsto \mathrm{vec}(\bfX\bfOmega^*)$ yields the desired result.
\end{proof}

Note that $k_0$ in \Cref{prop:ose gaussian trick forward} can be taken independently of the full dimension $q$, using for example rescaled Gaussian distribution, as detailed in \Cref{prop:ose with dim free gaussian}.
Hence, no additional dependency in the full dimension $q$ is introduced in the oblivious subspace embedding property of $\bfX \mapsto \bfX \bfOmega^*$, as detailed in \Cref{prop:ose gaussian trick} below.

\begin{proposition}
\label{prop:ose gaussian trick}
    Let $\varepsilon \in (0,1)$, $\delta \in (0,1)$ and $d\in\bbN\setminus\{0\}$.
    Let $\bfOmega$ be an $(\varepsilon', \delta', d)$ oblivious $\bbK^p \rightarrow \bbK^k$ subspace embedding, with 
    \[
        \varepsilon':= 0.8 \varepsilon,
        \quad 
        \delta' := \delta \lceil 8400 \varepsilon^{-2}(6.9 \eta_{\bbK} d + \log(8/\delta)) \rceil^{-1},
        \quad \eta_{\bbR}=1,
        \quad \eta_{\bbC}=1.
    \]
    Then, $\bfX \mapsto \bfX \bfOmega^*$ is an $(\varepsilon, \delta, d)$ oblivious $\bbK^{q\times p} \rightarrow \bbK^{q \times k}$ subspace embedding.
\end{proposition}
\begin{proof}
    Let $\bfOmega$ be an $(\varepsilon', \delta', d)$ oblivious subspace embedding with $\varepsilon ' := (1-\alpha) \varepsilon$, $\alpha \in (0,1)$, and $\delta'\in (0,1)$.
    Let $\varepsilon_0:=\alpha \varepsilon/4$, such that $\frac{1}{1-\varepsilon_0}(2\varepsilon_0 + \varepsilon' (1+\varepsilon_0)) \leq \varepsilon$, and let $\delta_0:=\delta/4$.
    From \Cref{prop:ose with dim free gaussian}, there exists a Gaussian $(\varepsilon_0, \delta_0, d)$ oblivious $\bbK^{p\times q} \rightarrow \bbK^{k_0\times q}$ subspace embedding with 
    \[
        k_0 := \lceil 10.5 \varepsilon_0^{-2}(6.9 \eta_{\bbK} d + \log(2/\delta_0)) \rceil
        = \lceil 10.5 (4/\alpha)^2\varepsilon^{-2}(6.9 \eta_{\bbK} d + \log(8/\delta)) \rceil.
    \]
    Now, from \Cref{prop:ose gaussian trick forward}, we obtain the desired result if $k_0 \delta' + 2\delta_0 \leq \delta$, or equivalently if
    \[
        \delta' 
        \leq \delta/2k_0
        = \delta \frac{1}{2\lceil 10.5 (4/\alpha)^2\varepsilon^{-2}(6.9 \eta_{\bbK} d + \log(8/\delta)) \rceil}
        \leq 
        \delta \frac{1}{\lceil 336 \alpha^{-2} \varepsilon^{-2}(6.9 \eta_{\bbK} d + \log(8/\delta)) \rceil}.
    \]
    Then, we would want to choose $\alpha$ such as to minimize the number of rows of $\bfOmega$.
    Although it depends on the class of embeddings considered, we can anticipate that $\bfOmega$ will have about $\calO(\varepsilon'^{-2}(d + \log(\delta')))$ rows.
    This drives us to choose $\alpha = 1/5$, and completes the proof.
\end{proof}


\section{Embedding for matrices}
\label{apdx:general ose matrix}

In this section we show how the oblivious embedding results on one side from \Cref{apdx:gaussian ose one side} and \Cref{apdx:general ose one side} yields an oblivious embedding property on $\bfUpsilon : \bfX \mapsto \bfGamma \mathrm{vec}(\bfOmega \bfX \bfSigma^*)$, as stated in \Cref{lem:merge three matrix ose} below.

\begin{lemma}
\label{lem:merge three matrix ose}
    Let $\bfSigma\in\bbK^{k_{\bfSigma} \times p}$ be an $(\varepsilon_{\bfSigma}, \delta_{\bfSigma}, d)$ oblivious $\bbK^{q\times p} \rightarrow \bbK^{q \times k_{\bfSigma}}$ subspace embedding.
    Let $\bfOmega\in\bbK^{k_{\bfOmega} \times q}$ be an $(\varepsilon_{\bfOmega}, \delta_{\bfOmega}, d)$ oblivious $\bbK^{q\times k_{\bfSigma}} \rightarrow \bbK^{k_{\bfOmega} \times k_{\bfSigma}}$ subspace embedding.
    Let $\bfGamma\in\bbK^{k_{\bfGamma} \times k_{\bfOmega} k_{\bfSigma}}$ be an $(\varepsilon_{\bfGamma}, \delta_{\bfGamma}, d)$ oblivious $\bbK^{k_{\bfOmega} k_{\bfSigma}} \rightarrow \bbK^{k_{\bfGamma}}$ subspace embedding.
    Then, the random linear map $\bfUpsilon : \bbK^{q\times p} \rightarrow \bbK^{k_{\bfGamma}}$ defined by
    \[
        \bfUpsilon(\bfX) 
        := \bfGamma \mathrm{vec}(\bfOmega \bfX \bfSigma^*)
    \]
    is an $(\varepsilon', \delta', d)$ oblivious $\bbK^{q\times p} \rightarrow \bbK^{k_{\bfGamma}}$ subspace embedding, with
    \[
        \varepsilon' := (1+\varepsilon_{\bfGamma})(1+\varepsilon_{\bfOmega})(1+\varepsilon_{\bfSigma}) - 1,
        \quad
        \delta' := \delta_{\bfGamma} + \delta_{\bfSigma} + \delta_{\bfOmega}.
    \]
\end{lemma}
\begin{proof}
    Let $V \subset \bbK^{q \times p}$ be a $d$-dimensional subspace.
    Using the embedding property of $\bfGamma$ and the independence of $\bfGamma$, $\bfSigma$ and $\bfOmega$, we obtain that
    \[
    \begin{aligned}
        \forall \bfX \in V,~ & \left\vert \|\bfX\|_F^2 - \| \bfUpsilon(\bfX)  \|_2^2 \right\vert
        = \left\vert \|\bfX\|_F^2 - \|\bfGamma \mathrm{vec}(\bfOmega\bfX\bfSigma^*) \|_2^2 \right\vert
        \\
        &\leq \left\vert \|\bfX\|_F^2 - \|\bfOmega\bfX\bfSigma^*\|_F^2 \right\vert
        +  \left\vert \|\bfOmega\bfX\bfSigma^*\|_F^2 - \|\bfGamma \mathrm{vec}(\bfOmega\bfX\bfSigma^*) \|_2^2 \right\vert
        \\
        & \leq \left\vert \|\bfX\|_F^2 - \|\bfOmega\bfX\bfSigma^*\|_F^2 \right\vert
        +  \varepsilon_{\bfGamma} \|\bfOmega\bfX\bfSigma^*\|_F^2
        \\
        & \leq (1+\varepsilon_{\bfGamma}) \left\vert \|\bfX\|_F^2 - \|\bfOmega\bfX\bfSigma^*\|_F^2 \right\vert + \varepsilon_{\bfGamma}\|\bfX\|_F^2
    \end{aligned}
    \]
    holds with probability at least $1-\delta_{\bfGamma}$.
    Then, using the embedding property of $\bfOmega$ and the independence of $\bfSigma$ and $\bfOmega$, we obtain that
    \[
        \begin{aligned}
            \forall \bfX \in V,~ 
            \left\vert 
                \|\bfX\|_F^2 - \|\bfOmega\bfX\bfSigma^*\|_F^2 
            \right\vert
            & \leq \left\vert 
                \|\bfX\|_F^2 - \|\bfX\bfSigma^*\|_F^2
            \right\vert
            + \left\vert 
                \|\bfX\bfSigma^*\|_F^2 - \|\bfOmega\bfX\bfSigma^*\|_F^2 
            \right\vert
            \\
            & \leq \left\vert 
                \|\bfX\|_F^2 - \|\bfX\bfSigma^*\|_F^2
            \right\vert
            + \varepsilon_{\bfOmega}\|\bfX\bfSigma^*\|_F^2
            \\
            & \leq (1+\varepsilon_{\bfOmega})\left\vert 
                \|\bfX\|_F^2 - \|\bfX\bfSigma^*\|_F^2
            \right\vert
            + \varepsilon_{\bfOmega}\|\bfX\|_F^2
        \end{aligned}
    \]
    holds with probability at least $1-\delta_{\bfOmega}$.
    Then, using the embedding property of $\bfSigma$, we obtain that
    \[
        \forall \bfX \in V,~ 
        \left\vert 
            \|\bfX\|_F^2 - \|\bfX\bfSigma^*\|_F^2 
        \right\vert
        \leq \varepsilon_{\bfSigma} \|\bfX\|_F^2
    \]
    holds with probability least $1-\delta_{\bfSigma}$.
    As a result, combining these last three results, and using the independence of $\bfGamma$, $\bfOmega$ and $\bfSigma$, we obtain that 
    \[
        \forall \bfX \in V,~ \big| \|\bfX\|_F^2 - \| \bfUpsilon(\bfX)  \|_2^2 \big|
        \leq 
        ((1+\varepsilon_{\bfGamma})(1+\varepsilon_{\bfOmega})(1+\varepsilon_{\bfSigma}) - 1) \|\bfX\|_F^2
        = \varepsilon' \|\bfX\|_F^2,
    \]
    with $\varepsilon' = (1+\varepsilon_{\bfGamma})(1+\varepsilon_{\bfOmega})(1+\varepsilon_{\bfSigma})-1$, holds with probability at least $1-\delta'$, with $\delta'=\delta_{\bfGamma} + \delta_{\bfSigma} + \delta_{\bfOmega}$.
    Finally, \Cref{lem:approx norm implies approx inner} yields the desired result.
\end{proof}

The above \Cref{lem:merge three matrix ose} can be for example used with various choices of accuracies and probabilities of failure.
A simple choice is $\varepsilon_{\bfSigma}=\varepsilon_{\bfOmega} = \varepsilon_{\bfGamma}$ and $\delta_{\bfSigma}=\delta_{\bfOmega}=\delta_{\bfGamma}$, which yields \Cref{prop:merge three matrix ose same prop}.
Note that it may be worth considering alternative choices for example when one wants to minimize the number of rows of $\bfSigma$ in order to reduce the number of matrix-vector products with $\bfX$.

\begin{proposition}
\label{prop:merge three matrix ose same prop}
    Let $\bfSigma\in\bbK^{k_{\bfSigma} \times p}$ be an $(\varepsilon/4, \delta/3, d)$ oblivious $\bbK^{q\times p} \rightarrow \bbK^{q \times k_{\bfSigma}}$ subspace embedding.
    Let $\bfOmega\in\bbK^{k_{\bfOmega} \times q}$ be an $(\varepsilon/4, \delta/3, d)$ oblivious $\bbK^{q\times k_{\bfSigma}} \rightarrow \bbK^{k_{\bfOmega} \times k_{\bfSigma}}$ subspace embedding.
    Let $\bfGamma\in\bbK^{k_{\bfGamma} \times k_{\bfOmega} k_{\bfSigma}}$ be an $(\varepsilon/4, \delta/3, d)$ oblivious $\bbK^{k_{\bfOmega} k_{\bfSigma}} \rightarrow \bbK^{k_{\bfGamma}}$ subspace embedding.
    Then, the random linear map $\bfUpsilon : \bbK^{q\times p} \rightarrow \bbK^{k_{\bfGamma}}$ defined by
    \[
        \bfUpsilon(\bfX) 
        := \bfGamma \mathrm{vec}(\bfOmega \bfX \bfSigma^*)
    \]
    is an $(\varepsilon, \delta, d)$ oblivious $\bbK^{q\times p} \rightarrow \bbK^{k_{\bfGamma}}$ subspace embedding.
\end{proposition}
\begin{proof}
    By \Cref{lem:merge three matrix ose}, $\bfUpsilon$ is an $(\varepsilon', \delta', d)$ oblivious $\bbK^{q\times p} \rightarrow \bbK^{k_{\bfGamma}}$ subspace embedding, with $\varepsilon' = (1+\varepsilon/4)^3 - 1$ and $\delta' = 3\delta/3 = \delta$.
    Now, by convexity of $h(t) := (1+t/4)^3-1$ on $(0,1)$ we have $\varepsilon' = h(\varepsilon) \leq (1-\varepsilon) h(0) + \varepsilon h(1) = \frac{61}{64}\varepsilon \leq \varepsilon$, which yields the desired result.
\end{proof}

\bibliographystyle{plain}  
\bibliography{main}

\end{document}